\documentclass[11pt]{article}
\usepackage{amssymb}
\usepackage{amsthm}
\usepackage{amsmath}
\usepackage{geometry}
\usepackage{enumitem}
\usepackage{url}
\usepackage{mathrsfs}
\usepackage{dsfont}
\usepackage[hidelinks]{hyperref}

\numberwithin{equation}{section}

\newtheorem{thm}{Theorem}[section]
\newtheorem{prop}[thm]{Proposition}
\newtheorem{lem}[thm]{Lemma}
\newtheorem{cor}[thm]{Corollary}

\theoremstyle{definition}

\newcommand{\diff}{\mathop{}\!\mathrm{d}}
\newcommand{\scp}[2]{\langle #1, #2\rangle}
\newcommand{\e}{\mathrm{e}}
\newcommand{\I}{\mathrm{i}}

\newcommand{\angles}[1]{\langle #1 \rangle}

\DeclareMathOperator*{\slim}{s-lim}

\title{Resolvent Boundary Values and Higher-Order Regularity Estimates for Unitary and Self-Adjoint Operators}

\author{Janik Kruse\textsuperscript{1,2} \vspace{0.2cm} \\
	\small \textsuperscript{1}Technische Universität Braunschweig\\
	\small Carl-Friedrich-Gauß-Fakultät, Institut für Analysis und Algebra \\
	\small Universitätsplatz 2,	38106 Braunschweig, Germany \vspace{0.2cm} \\
	\textsuperscript{2}\small Nagoya University \\ 
	\small Institute of Liberal Arts and Sciences \\
	\small Furo-cho, Chikusa-ku, Nagoya 464-8601, Japan}
	
\date{November 26, 2026}

\begin{document}

\maketitle

\begin{abstract}
	We study the boundary values of the resolvent of a unitary operator using Mourre's conjugate operator method. Under the regularity assumption $U\in\mathcal C^{1,1}(A)$, where $A$ is	the self-adjoint conjugate operator, we establish the existence	and weak$^*$ continuity of these boundary values on the Mourre set $\mu^A(U)$ as operators from $\mathcal H_{1/2,1}$ to $\mathcal H_{-1/2,\infty}$. If	$U\in\mathcal C^{s+1/2}(A)$ for $s>1/2$, the boundary-value functions are locally of Lipschitz--Zygmund class $\Lambda^{s-1/2}$ as operators from $\mathcal H_s$ to	$\mathcal H_{-s}$. The proofs use the method of differential inequalities adapted to unitary operators. The results are	optimal within the regularity and weighted-space scales	considered. We also recover the corresponding results for self-adjoint operators, but without imposing a spectral-gap assumption as in previous results.
\end{abstract}

\section{Introduction}

The analysis of the existence and regularity of resolvent boundary values is an important part of spectral analysis. For self-adjoint operators, the boundary values have been extensively studied within Mourre's conjugate operator theory \cite{mourre1983,amrein1996,boutetDeMonvel1996,sahbani1997, gerard2008}. In this framework, suitable regularity of a self-adjoint operator $H$ with respect to a self-adjoint conjugate operator~$A$, together with a positive commutator estimate, yields estimates and boundary regularity for the resolvent of~$H$ near its spectrum. The purpose of the present paper is to develop boundary-value and higher-order regularity results for resolvents of unitary operators, following the method of differential inequalities used in \cite{boutetDeMonvel1996} for self-adjoint operators.

Let $U$ be a unitary operator on a Hilbert space $\mathcal H$. We study the operator-valued function $z\mapsto (1-zU^*)^{-1}$, which is holomorphic inside and outside the unit circle $\mathbb T=\{z\in\mathbb C\mid |z|=1\}$. However, as $z$ approaches a point of the spectrum from either side of the unit circle, the resolvent diverges in the operator norm topology on $\mathcal H$. Boundary values of the resolvent must therefore be analysed between suitable weighted spaces. The conjugate operator method provides a natural scale of Sobolev and Besov spaces associated with the conjugate operator $A$, on which the resolvent near the spectrum can be studied.

Our first result (Theorem~\ref{thm:LAP}) establishes the existence and weak$^*$ continuity of the boundary values under the regularity assumption $U\in\mathcal C^{1,1}(A)$. The limits are taken from either side of the unit circle on the region where a strict Mourre estimate holds. The boundary values define bounded operators from the Besov space $\mathcal{H}_{1/2,1}$ associated with $A$ to its dual $\mathcal{H}_{-1/2,\infty}$. Our second result (Theorem~\ref{thm:HigherRegularity}) describes the regularity of these boundary values. Under the stronger assumption $U\in\mathcal{C}^{s+1/2}(A)$, for $s>1/2$, the boundary values belong to the Lipschitz--Zygmund class~$\Lambda^{s-1/2}$, as operators from the Sobolev space $\mathcal{H}_s$ to its dual $\mathcal{H}_{-s}$. We refer to Section~\ref{ssec:Notation} below for the definitions of these spaces and regularity classes.

The proof follows the method of differential inequalities as presented in \cite{boutetDeMonvel1996}. First, we regularise the unitary operator~$U$ along the conjugation group generated by $A$. The regularisation is chosen such that, in leading order, the defect of the regularised operator of being unitary is proportional to the commutator appearing in the Mourre estimate. This allows us to prove the existence of the regularised resolvent on the unit circle. Then, we prove that differentiation with respect to the regularisation parameter produces commutator identities with error terms controlled by the regularity of $U$. These identities give the existence of the resolvent boundary values and the estimates for higher-order derivatives needed to obtain Lipschitz--Zygmund regularity.

In previous works, some aspects of the Mourre theory for unitary operators have already been investigated. For example, in \cite{fernandez2013}, it was observed that much of the existing self-adjoint theory can be transferred to unitary operators by means of the Cayley transform. Moreover, in \cite{astaburuaga2015}, the existence of the boundary values is stated under the regularity assumption $U\in\mathcal C^{1,1}(A)$, and higher-order regularity results are discussed in the preprint~\cite{astaburuaga2013}. Here, we give an independent proof of the existence theorem and establish higher-order regularity results which are stronger than the ones stated in \cite{astaburuaga2013}.

Our approach also allows us to recover the corresponding theorems for self-adjoint operators, but without imposing a spectral-gap assumption as in \cite[Theorem~7.4.1]{amrein1996} and \cite[Theorem on p.~12]{boutetDeMonvel1996}. Results for self-adjoint operators without a spectral-gap assumption were previously obtained under stronger regularity assumptions. The best preceding ours is \cite{sahbani1997}, where the class $C^{1+0}(A)$ was used for the existence of the boundary values. Because this class is strictly contained in $\mathcal C^{1,1}(A)$, our results close a gap in the self-adjoint theory. Moreover, they strengthen several results which rely on the existence and regularity of the boundary values, including \cite[Corollary~3.6]{kruse2024} and \cite[Section~2]{fernandez2013}. We remark that another proposed strategy for removing the spectral-gap assumption, which is mentioned in \cite[p.~248]{amrein1996}, does not work, as we illustrate in Appendix~\ref{appx:Counterexample}.

Removing the spectral-gap assumption is relevant in applications where the spectrum consists of the entire real line in the self-adjoint setting, or the entire unit circle in the unitary setting. Typical examples include Schrödinger operators with electric potentials \cite{jensen1987} and unitary operators describing the evolution of quantum walks \cite{richard2018}. The higher-order regularity result also makes it possible to study the regularity and mapping properties of the scattering operator, using the stationary representation formula of \cite{tiedra2020}.

\subsection{Notation}
\label{ssec:Notation}

We use notation and terminology that are standard in Mourre theory. We therefore recall only the main definitions and refer to the literature for further details \cite{amrein1996,boutetDeMonvel1996}.

Throughout the paper, $\mathcal H$ denotes a complex separable Hilbert space, $U$ a unitary operator on~$\mathcal H$, and $A$ a self-adjoint operator on $\mathcal H$. For Banach spaces $\mathcal X$ and $\mathcal Y$, we write $B(\mathcal X,\mathcal Y)$ for the space of bounded linear operators from $\mathcal X$ to $\mathcal Y$, and set $B(\mathcal H)=B(\mathcal H,\mathcal H)$. We identify the Hilbert space $\mathcal{H}$ with its dual space $\mathcal{H}^*$. 

\subsubsection{Sobolev and Besov Spaces}
\label{ssec:SobolevBesovSpaces}

We write $\angles{A} = (1+A^2)^{1/2}$ and denote by $\{\mathcal H_s\}_{s\in\mathbb R}$ the Sobolev scale associated with $A$. For $s\geq0$, the space $\mathcal H_s$ is the domain of $\langle A\rangle^s$, equipped with the norm $\|f\|_s = \|\langle A\rangle^s f\|$. For $s<0$, the space $\mathcal H_s$ is the completion of $\mathcal H$ with respect to the same norm. We also use the Besov scale $\{\mathcal H_{s,p}\mid s\in\mathbb R,\ p\in[1,\infty] \}$ obtained by real interpolation of the Sobolev scale. More precisely, if $s_0<s<s_1$, we set
\begin{align}
	\mathcal H_{s,p} = (\mathcal H_{s_1},\mathcal H_{s_0})_{\theta,p} \, , \qquad \theta = \frac{s_1-s}{s_1-s_0},
\end{align}
where $(\cdot,\cdot)_{\theta,p}$ denotes the real interpolation functor defined by the $K$-method. These spaces are independent, up to equivalence of norms, of the choice of $s_0$ and $s_1$ \cite[Section~3.4.1]{amrein1996}.

The spaces relevant to our main results are $\mathcal H_{1/2,1}$ and $\mathcal H_s$, with $s>1/2$, together with their duals $\mathcal{H}_{-1/2,\infty}$ and $\mathcal{H}_{-s}$, respectively. We also use the continuous embedding $\mathcal{H}_s = \mathcal{H}_{s,2} \subset \mathcal{H}_{1/2,1}$, for $s>1/2$.

\subsubsection{Operator Regularity}
\label{ssec:OperatorRegularity}

Let $W_\tau=\e^{\I\tau A}$, for $\tau\in\mathbb R$, be the unitary group generated by $A$. Conjugation by $W_\tau$ defines a group of automorphisms on $B(\mathcal H)$ given by $\mathcal W_\tau[T] = W_{-\tau}TW_\tau$. For every $T\in B(\mathcal H)$, the map $\tau\mapsto\mathcal W_\tau[T]$ is continuous in the strong operator topology.

For $k\in\mathbb N$, we say that $T\in C^k(A)$ if the map $\tau\mapsto \mathcal{W}_\tau[T]$ is of class $C^k$ in the strong operator topology, and that $T\in C^k_{\mathrm u}(A)$ if it is of class $C^k$ in the operator norm topology. We equip the space $C_\mathrm{u}^k(A)$ with the norm $\|T\|_{C_\mathrm{u}^k(A)} = \sum_{j=0}^{k} \|\mathcal A^j[T]\|$.

The commutator $TA-AT$ is initially understood as a sesquilinear form on $D(A)$. For $T\in C^1(A)$, this form extends to a bounded operator on $\mathcal H$, which we denote by $\mathcal A[T] =[T,A]$ \cite[Proposition~5.1.2]{amrein1996}. With this convention, we have, for $T\in C^1(A)$,
\begin{align}
	\frac{\diff}{\diff\tau}\mathcal W_\tau[T] = \I\mathcal W_\tau \mathcal A[T],
\end{align}
where the derivative is taken in the strong operator topology. 

To describe operator classes of fractional regularity, we interpolate the Banach spaces $C_\mathrm{u}^k(A)$. For $s>0$ and $p\in[1,\infty]$, we choose an integer $\ell>s$ and define
\begin{align}
	\mathcal C^{s,p}(A)	= (C_\mathrm{u}^\ell(A), C_\mathrm{u}^0(A))_{\theta,p} \, , \qquad \theta = 1-\frac{s}{\ell}.
\end{align}
These spaces are independent, up to equivalence of norms, of the integer $\ell>s$ \cite[Section~5.2.1]{amrein1996}, and admit equivalent characterisations by finite differences of $\mathcal W_\tau$ \cite[Section~2.1]{boutetDeMonvel1996}. We abbreviate $\mathcal C^s(A)=\mathcal C^{s,\infty}(A)$. Note that the classes $\mathcal C^k(A)$ should be distinguished from the classes $C^k(A)$ at integer orders. We also use the continuous embeddings $\mathcal C^s(A)\subset\mathcal C^{1,1}(A)\subset C_\mathrm{u}^1(A)$, for $s>1$.

\subsubsection{Lipschitz--Zygmund Class}
\label{ssec:LPClass}

Let $\mathcal{X}$ be a Banach space. For a function $f:\mathbb{T}\to\mathcal{X}$ and $\theta\in\mathbb{R}$, we define the finite difference operator
\begin{align} \label{eq:FiniteDifference}
	\Delta_{\theta} f(\zeta) = f(\e^{\I\theta} \zeta) - f(\zeta),
\end{align}
and denote its $m$-fold iteration, for $m\in\mathbb{N}$, by $\Delta_\theta^m$. The tangential derivative of $f$ is defined by
\begin{align}
	D_\zeta f(\zeta) = \lim_{\theta\to0} \frac{\Delta_\theta f(\zeta)}{\theta},
\end{align}
whenever this limit exists in $\mathcal{X}$. Note that $D_\zeta \zeta = \I \zeta$.

For $s>0$, we say that a continuous function $f$ belongs to the Lipschitz--Zygmund class~$\Lambda^s$ if there exist an integer $m>s$ and a constant $c<\infty$ such that, for all $\zeta\in \mathbb{T}$ and $\theta\in\mathbb{R}$,
\begin{align}
	\|\Delta^m_\theta f(\zeta)\|_{\mathcal{X}} \leq c |\theta|^s.
\end{align}
If the function $f$ is defined on an open subset $I\subset\mathbb{T}$, we say that $f$ is locally of class $\Lambda^s$ if, for all $\eta\in C_0^\infty(I)$, the function $\eta f$, extended by zero outside $I$, is of class $\Lambda^s$.

\subsection{The Mourre Estimate}

Assume that $U\in C^1(A)$. Then, $[A,U]U^*$ and $U^*[A,U]$ are bounded and self-adjoint operators on $\mathcal H$. We say that a strict Mourre estimate holds at $\zeta\in\mathbb T$ if there exist a real-valued function $\eta\in C^\infty(\mathbb T)$, with $\eta(\zeta)\neq 0$, and a constant $a>0$ such that
\begin{subequations} \label{eq:MourreEstimate}
\begin{align}
	\eta(U)[A,U]U^*\eta(U) \geq a\,\eta(U)^2. 
\end{align}
Conjugating by $U$, and using that $U$ commutes with $\eta(U)$, gives the equivalent estimate
\begin{align}
	\eta(U)U^*[A,U]\eta(U) \geq a\,\eta(U)^2.
\end{align}
\end{subequations}
The Mourre set $\mu^A(U)$ of $U$ with respect to $A$ is defined as the set of all $\zeta\in\mathbb{T}$ for which a Mourre estimate holds at $\zeta$. It is easily seen that $\mu^A(U)$ is an open subset of $\mathbb{T}$. Indeed, if $\zeta_0\in \mu^A(U)$ and \eqref{eq:MourreEstimate} holds with $\eta(\zeta_0)\neq 0$, then the same function $\eta$ and constant $a$ establish the estimate at every point of the open neighbourhood $\{\zeta\in\mathbb T\mid\eta(\zeta)\ne0\}$ of $\zeta_0$.

More generally, Mourre estimates may include a compact remainder. At a point which is not an eigenvalue, such an estimate yields a Mourre estimate without a compact remainder after shrinking the support of~$\eta$ \cite[Lemma~7.2.12]{amrein1996}.

There are plenty of examples in the literature where Mourre estimates for unitary operators have been established. An elementary example is the bilateral shift on $\ell^2(\mathbb Z)$, for which the number operator provides a conjugate operator \cite[Section~3.1]{fernandez2013}. Indeed, let $\{e_n\}_{n\in\mathbb{Z}}$ be the standard basis of $\ell^2(\mathbb{Z})$. Then, with $Ue_n=e_{n+1}$ and $Ae_n=ne_n$, we get $[A,U]U^*=1$, and hence $\mu^A(U)=\mathbb T$. Further examples of Mourre estimates for unitary operators include Koopman operators \cite{tiedra2015}, Floquet operators \cite{asch2016}, and quantum walks with position-dependent coins \cite{richard2018}.

\subsection{Main Results}

Our first result establishes the existence and weak$^*$ continuity of the boundary values of the resolvent $(1-zU^*)^{-1}$ on the Mourre set $\mu^A(U)$ under the regularity assumption $U\in\mathcal C^{1,1}(A)$. The limits are taken from either side of the unit circle, through the domains $\mathbb D_+=\{z\in\mathbb C\mid |z|<1\}$ and $\mathbb D_-=\{z\in\mathbb C\mid |z|>1\}$, respectively. The limits can be taken arbitrarily within $\mathbb{D}_\pm$, and the two boundary values need not coincide.

\begin{thm} \label{thm:LAP}
	Let $U\in\mathcal C^{1,1}(A)$ be a unitary operator on $\mathcal H$. Then, the holomorphic functions
	\begin{align} \label{eq:LAP}
		\mathbb D_\pm\ni z \mapsto (1-zU^*)^{-1} \in B(\mathcal H_{1/2,1},\mathcal H_{-1/2,\infty})
	\end{align}
	extend to weak$^*$ continuous functions on $\mathbb D_\pm\cup\mu^A(U)$.
\end{thm}

Here, weak$^*$ continuity means that, for every $f,g\in\mathcal H_{1/2,1}$, the holomorphic functions $\mathbb{D}_\pm \ni z\mapsto\scp{g}{(1-zU^*)^{-1}f}$ extend to continuous functions on $\mathbb D_\pm\cup\mu^A(U)$. We denote the boundary values of \eqref{eq:LAP} by $B_\pm(\zeta)$, for $\zeta \in \mu^A(U)$. These limits define bounded operators from $\mathcal H_{1/2,1}$ to $\mathcal H_{-1/2,\infty}$.

The theorem is optimal with respect to the regularity assumption and the spaces involved. In general, the condition
$U\in\mathcal C^{1,1}(A)$ cannot be weakened to $U\in\mathcal C^{\alpha,p}(A)$ with $\alpha<1$, or with $\alpha=1$ and $p>1$, and the domain $\mathcal H_{1/2,1}$ cannot be replaced by $\mathcal H_{s,p}$ with $s<1/2$, or with $s=1/2$ and $p>1$. Nor can weak$^*$ continuity be strengthened to strong operator continuity. These optimality statements are known in the self-adjoint setting \cite[p.~12]{boutetDeMonvel1996} and carry over to unitary operators by the results of Section~\ref{ssec:SelfAdjointOperators}.

For $s>1/2$, the embedding $\mathcal{H}_s \subset \mathcal{H}_{1/2,1}$ allows us to regard the boundary values $B_\pm(\zeta)$ as bounded operators from $\mathcal H_s$ to $\mathcal H_{-s}$. Our second result describes their regularity in $\zeta$ under stronger regularity assumptions on~$U$.

\begin{thm} \label{thm:HigherRegularity}
	Let $s>1/2$, and let $U\in\mathcal C^{s+1/2}(A)$ be a unitary operator on $\mathcal{H}$. Then, the boundary-value functions
	\begin{align}
		\mu^A(U) \ni \zeta \mapsto B_\pm(\zeta) \in B(\mathcal{H}_s,\mathcal{H}_{-s})
	\end{align}
	are locally of class $\Lambda^{s-1/2}$.
\end{thm}

Here, $\Lambda^{s-1/2}$ denotes the Lipschitz--Zygmund class on $\mathbb{T}$ defined in Section~\ref{ssec:LPClass}. Again, as is known from the self-adjoint case, the theorem is optimal within the scales $\mathcal C^t(A)$, $\mathcal H_s$, and $\Lambda^\alpha$ \cite[p.~12]{boutetDeMonvel1996}.

\subsection{Self-Adjoint Operators}
\label{ssec:SelfAdjointOperators}

We explain how our main results for unitary operators yield the corresponding results for self-adjoint operators. The connection is provided by the Cayley transform, which was first observed in \cite{fernandez2013}. This allows us to recover the results of \cite[Theorem~7.4.1]{amrein1996} and \cite[Theorem on p.~12]{boutetDeMonvel1996}, but without imposing a spectral-gap assumption.

The Cayley transform $\kappa:\mathbb C\backslash\{-\I\} \to \mathbb{C}$ is defined by
\begin{align}
	\kappa(z) = \frac{z-\I}{z+\I}.
\end{align}
For $z$ in the upper half-plane $\mathbb C_+=\{z\in\mathbb C\mid \Im z>0\}$, the point $\kappa(z)$ belongs to $\mathbb D_+$, and for $z$ in the lower half-plane $\mathbb C_-=\{z\in\mathbb C\mid \Im z<0\}$, the point $\kappa(z)$ belongs to $\mathbb D_-$, apart from the pole of $\kappa$ at $z=-\I$.

Let $H$ be a self-adjoint operator on $\mathcal H$. We recall that the regularity classes $\mathcal C^{s,p}(A)$ for unbounded self-adjoint operators are understood in the sense that $H$ is of class $\mathcal C^{s,p}(A)$ if its resolvent $(H+\I)^{-1}$ is of class $\mathcal C^{s,p}(A)$ \cite[Definition~6.2.2]{amrein1996}. We define the Cayley transform of $H$ by
\begin{align}
	U = (H-\I)(H+\I)^{-1}.
\end{align}
The operator $U$ is unitary and satisfies $U=1-2\I(H+\I)^{-1}$. Thus, the regularity of $H$ with respect to $A$ is equivalent to the corresponding regularity of its Cayley transform $U$, that is, $H\in\mathcal C^{s,p}(A)$ if and only if $U\in\mathcal C^{s,p}(A)$.

We recall the Mourre set of $H$ with respect to $A$. If $H\in C^1(A)$, the commutator $[H,\I A] = \I(HA-AH)$, initially defined as a sesquilinear form on $D(A)\cap D(H)$, extends continuously to a sesquilinear form on $D(H)$ \cite[Theorem~6.2.10]{amrein1996}. We denote by $\mu^A(H)$ the set of all $\lambda\in\mathbb R$ for which there exist a real-valued function $\chi\in C_0^\infty(\mathbb R)$, with $\chi(\lambda)\ne0$, and a constant $a>0$ such that
\begin{align}
	\chi(H)[H,\I A]\chi(H)\geq a\,\chi(H)^2.
\end{align}
The commutators associated with $H$ and $U$ are related by the formula
\begin{align}
	(H+\I)^{-1}[H,\I A](H-\I)^{-1} = \frac{1}{2}[U,A]U^*.
\end{align}
Consequently, the Mourre estimates for $H$ and for its Cayley transform $U$ are equivalent, in the sense that the Mourre sets satisfy the identity
\begin{align}
	\kappa(\mu^A(H)) = \mu^{-A}(U)\backslash\{1\}.
\end{align}
Moreover, for $z\in\mathbb{C}\backslash(\mathbb{R} \cup \{\pm\I\})$, the resolvents of $H$ and $U$ are related by
\begin{align} \label{eq:ResolventIdentity}
	(H-z)^{-1} =\frac{2\I}{z^2+1}(1-\kappa(z)U^*)^{-1} - \frac{1}{z-\I}.
\end{align}
Together with Theorem~\ref{thm:LAP}, this identity yields the existence and weak$^*$ continuity of the boundary values of $(H-z)^{-1}$ as $z$ approaches the Mourre set of $H$ from the upper or lower half-plane. Note that the conjugate operator $-A$ defines the same Sobolev and Besov spaces as $A$.

\begin{thm}
	Let $H\in\mathcal C^{1,1}(A)$ be a self-adjoint operator on $\mathcal H$. Then, the holomorphic functions
	\begin{align} \label{eq:LAPSA}
		\mathbb C_\pm\ni z\mapsto (H-z)^{-1} \in B(\mathcal H_{1/2,1},\mathcal H_{-1/2,\infty})
	\end{align}
	extend to weak$^*$ continuous functions on $\mathbb C_\pm\cup\mu^A(H)$.
\end{thm}

We denote the boundary values of \eqref{eq:LAPSA} by $R_\pm(\lambda)$, for $\lambda\in\mu^A(H)$. According to the resolvent identity \eqref{eq:ResolventIdentity}, they are related to the boundary values $B_\pm$ of the Cayley transform of $H$ by
\begin{align}
	R_\pm(\lambda) = \frac{2\I}{\lambda^2+1} B_\pm(\kappa(\lambda))
	-\frac{1}{\lambda-\I}.
\end{align}
For functions on~\(\mathbb R\), we use the Lipschitz–Zygmund classes defined by translation differences in \cite[Section~1.2]{boutetDeMonvel1996}. These classes are compatible with the classes defined in Section~\ref{ssec:LPClass} in the sense that the Cayley transform maps functions on $\mathbb{R}$ which are locally of Lipschitz--Zygmund class to functions on $\mathbb{T}\backslash\{1\}$ which are locally of Lipschitz--Zygmund class, and vice-versa. Thus, Theorem~\ref{thm:HigherRegularity} yields the following higher-order regularity result for the resolvent boundary values of self-adjoint operators.

\begin{thm}\label{thm:selfadjoint-regularity}
	Let $s>1/2$, and let $H\in \mathcal C^{s+1/2}(A)$ be a self-adjoint operator on $\mathcal{H}$. Then, the boundary-value functions
	\begin{align}
		\mu^A(H)\ni\lambda \mapsto R_\pm(\lambda) \in B(\mathcal H_s,\mathcal H_{-s})
	\end{align}
	are locally of class $\Lambda^{s-1/2}$.
\end{thm}

\section{Regularisation of $U$ with Respect to $A$}
\label{sec:Regularisation}

In this section, we construct a family of bounded operators $U_\varepsilon$ that approximates $U$ and is smooth with respect to $A$ for every $\varepsilon\ne0$. The regularisation is chosen such that its derivative with respect to the regularisation parameter $\varepsilon$ agrees with the commutator $\mathcal A[U_\varepsilon]$, up to an error whose behaviour near $\varepsilon=0$ characterises the regularity of $U$. This construction is taken from \cite[Sections~4.1--4.3]{boutetDeMonvel1996}.

Choose a real-valued function $\varphi \in C_0^\infty(\mathbb R)$ such that $\varphi=1$ on a neighbourhood of zero, and let $\xi_0\in C_0^\infty(\mathbb{R})$ be the function $\xi_0(x) = \e^x\varphi(x)$. For $\varepsilon\in\mathbb R$, we define the operators
\begin{align} \label{eq:Uepsilon}
	U_\varepsilon =	\xi_0(\varepsilon\mathcal A)[U]	
	= \frac{1}{\sqrt{2\pi}} \int_{\mathbb R} \e^{-\I\varepsilon\tau A} U \e^{\I\varepsilon\tau A} \widehat{\xi_0}(\tau)\diff\tau,
\end{align}
where the Fourier transform of $\xi_0$ is
\begin{align}
	\widehat{\xi_0}(\tau) = \frac1{\sqrt{2\pi}} \int_{\mathbb R} \e^{-\I\tau x}\xi_0(x) \diff x.
\end{align}
The integral in \eqref{eq:Uepsilon} is understood in the strong operator sense. From $\xi_0(0)=1$, it follows that $\slim_{\varepsilon\to 0} U_\varepsilon = U_0 = U$. For $\varepsilon\neq 0$, the operator $U_\varepsilon$ is not necessarily unitary, but $(U_\varepsilon)^* = (U^*)_{-\varepsilon}$. Moreover, $U_\varepsilon$ is uniformly bounded in $\varepsilon$. Indeed, its operator norm is bounded by
\begin{align}
	\|U_\varepsilon\| \leq \frac{1}{\sqrt{2\pi}} \|\widehat{\xi_0}\|_{L^1}\|U\|.
\end{align}

\subsection{Regularity of $U_\varepsilon$}

It is easy to verify that $\varepsilon\mapsto U_\varepsilon$ is of class $C^\infty$ on $\mathbb{R}\backslash\{0\}$, and that $U_\varepsilon\in C^\infty(A)$ for all $\varepsilon\neq 0$. Indeed, for $k\in\mathbb{N}$, let $\eta_k \in C_0^\infty(\mathbb{R})$ be the function $\eta_k(x) = x^k \xi_0^{(k)}(x)$. Then, for all $\varepsilon\neq0$, we have
\begin{align} \label{eq:DerivativesUepsilon}
	U_\varepsilon^{(k)} = \partial_\varepsilon^k U_\varepsilon = \mathcal{A}^k \xi_0^{(k)}(\varepsilon\mathcal{A})[U] = \varepsilon^{-k} \eta_k(\varepsilon\mathcal{A})[U].
\end{align}
Moreover, if $U\in C^k(A)$, then $\varepsilon\mapsto U_\varepsilon$ is of class $C^k$ on $\mathbb{R}$ in the strong operator topology. In this case, $\mathcal{A}^k$ can be applied to $U$, giving $U_\varepsilon^{(k)} = \xi_0^{(k)}(\varepsilon\mathcal{A}) \mathcal{A}^k[U]$. Moreover, $U_0^{(k)} = \mathcal{A}^k[U]$ because $\xi_0^{(k)}(0)=1$ for all $k\in\mathbb{N}$. 

\begin{lem} \label{lem:Uepsilon}
	Assume that $U\in C^1(A)$. Then, for every integer $j\geq 1$,
	\begin{align}
		\|U_\varepsilon-U\| &\leq c_0 |\varepsilon|,\\
		\|U_\varepsilon^{(j)}\|	&\leq c_j|\varepsilon|^{1-j}, \qquad \varepsilon\neq 0.
	\end{align}
\end{lem}

\begin{proof}
	From $U\in C^1(A)$, it follows that $\|U_\varepsilon'\| \leq (2\pi)^{-\frac{1}{2}}\|\widehat{\xi'_0}\|_{L^1} \|\mathcal{A}[U]\|$; hence,
	\begin{align} \label{eq:UepsilonRateOfConvergence}
		\| U_\varepsilon - U\| = \left\| \int_0^\varepsilon U_\tau' \diff \tau \right\| \leq c_0 |\varepsilon|.
	\end{align}
	Moreover, for $j\geq 1$, we have
	\begin{align}
		U_\varepsilon^{(j)} = \mathcal{A}^{j-1} \xi_0^{(j)}(\varepsilon\mathcal{A})\mathcal{A}[U] = \varepsilon^{1-j} (\varepsilon\mathcal{A})^{j-1} \xi_0^{(j)}(\varepsilon\mathcal{A}) \mathcal{A}[U],
	\end{align}
	where the function $x\mapsto x^{j-1}\xi_0^{(j)}(x)$ is of class $C_0^\infty(\mathbb{R})$ if $j\geq 1$.
\end{proof}

\subsection{First-Order Expansion and Unitarity Defects of $U_\varepsilon$}

Let us now assume that $U\in C^1_\mathrm{u}(A)$. Then, the function $\varepsilon\mapsto U_\varepsilon$ is of class $C^1$ in norm on $\mathbb{R}$. Because $\mathcal{A}[U]$ is of class $C_\mathrm{u}^0(A)$ and $U_\varepsilon' = \xi_0'(\varepsilon\mathcal A)\mathcal A[U]$, we have
\begin{align}
	\lim_{\varepsilon\to 0} \| U_\varepsilon' - \mathcal{A}[U] \| = 0.
\end{align}
Consequently, $U_\varepsilon = U + \varepsilon\mathcal A[U] + o(\varepsilon)$ in operator norm, and
\begin{subequations}\label{eq:ConvergenceToCommutator}
	\begin{align}
		\lim_{\varepsilon\to 0} \| \varepsilon^{-1} (\mathds{1} - U_{\varepsilon} (U_{\varepsilon})^*) + 2\mathcal{A}[U] U^* \| &= 0, \label{eq:ConvergenceToCommutator1}\\
		\lim_{\varepsilon\to 0} \| \varepsilon^{-1} (\mathds{1} - (U_{\varepsilon})^* U_{\varepsilon}) + 2U^*\mathcal{A}[U] \| &= 0. \label{eq:ConvergenceToCommutator2}
	\end{align}
\end{subequations}
The above results give an estimate on the \textit{unitarity defects} of $U_\varepsilon$. In leading order in $\varepsilon$, they are given by $2\varepsilon[A,U]U^*$ and $2\varepsilon U^*[A,U]$, which are proportional to the commutators appearing in the Mourre estimate.

\subsection{A Multiplier Characterisation of Regularity Classes}

Let $\xi_1\in C_0^\infty(\mathbb{R})$ be the function $\xi_1(x) = \e^x x \varphi'(x)$. Using \eqref{eq:DerivativesUepsilon}, we get that, for $\varepsilon\neq 0$,
\begin{align}
	U_\varepsilon' &= \mathcal A[U_\varepsilon]-K_\varepsilon, \label{eq:UepsilonPrime}\\
	K_\varepsilon &=-\varepsilon^{-1}\xi_1(\varepsilon\mathcal A)[U]. \label{eq:Kepsilon}
\end{align}
Thus, the operator $K_\varepsilon$ describes the difference between the commutator $\mathcal A[U_\varepsilon]$ and the derivative of the regularised operator $U_\varepsilon$. The map $\varepsilon\mapsto K_\varepsilon$ is of class
$C^\infty$ on $\mathbb R\backslash\{0\}$. The following proposition characterises the regularity of $U$ through the behaviour of $K_\varepsilon$ as $\varepsilon\to 0$. It also establishes corresponding estimates for all derivatives of $K_\varepsilon$, which will be useful later.

\begin{prop} \label{prop:Regularity}
	Let $s>0$, $p\in[1,\infty]$, and let $j\ge0$ be an integer. Then, $U\in\mathcal C^{s,p}(A)$ if and only if, for either choice of sign,
	\begin{align} \label{eq:MultiplierCharacterisation}
		\begin{cases}
			\displaystyle \int_0^1 \left(\varepsilon^{j+1-s}\|K_{\pm\varepsilon}^{(j)}\|\right)^p \frac{\diff\varepsilon}{\varepsilon}<\infty,
			& 1\le p<\infty,\\[1em]
			\displaystyle \sup_{0<\varepsilon<1} \varepsilon^{j+1-s}\|K_{\pm\varepsilon}^{(j)}\|<\infty,
			& p=\infty.
		\end{cases}
	\end{align}
\end{prop}

\begin{proof}	
	We use the multiplier characterisation of the classes $\mathcal C^{s,p}(A)$ given in \cite[Corollary~2.3]{boutetDeMonvel1996}. In the form needed here, it states that, if $\psi\in C_0^\infty(\mathbb{R}\backslash\{0\})$ is not identically zero on $(0,\infty)$ and on $(-\infty,0)$, then, for $T\in B(\mathcal{H})$, $s>0$, and $p\in [1,\infty)$,
	\begin{subequations} \label{eq:RegularityMultiplierCharacterisation}
		\begin{align}
			T\in\mathcal C^{s,p}(A)	\quad\Leftrightarrow\quad
			\int_0^1 \left(\varepsilon^{-s}\|\psi(\varepsilon\mathcal A)[T]\|\right)^p
			\frac{\diff\varepsilon}{\varepsilon}<\infty,
		\end{align}
		and, for $p=\infty$, 
		\begin{align}
			T\in\mathcal C^{s,\infty}(A) \quad\Leftrightarrow\quad
			\sup_{0<\varepsilon<1}	\varepsilon^{-s}\|\psi(\varepsilon\mathcal A)[T]\|<\infty.
		\end{align}	
	\end{subequations}
	We set $\psi_0=\xi_1$ and define $\psi_{j+1}(x)=x\psi_j'(x)-(j+1)\psi_j(x)$. Each $\psi_j$ belongs to $C_0^\infty(\mathbb R\backslash\{0\})$, and is not identically zero on $(0,\infty)$ and on $(-\infty,0)$. Moreover, for $\varepsilon>0$, it holds that
	\begin{align}
		\varepsilon^{j+1-s} \|K_{\pm\varepsilon}^{(j)}\| = -\varepsilon^{-s} \|\psi_j(\pm\varepsilon\mathcal A)[U]\|.
	\end{align}
	Thus, applying the multiplier characterisation \eqref{eq:RegularityMultiplierCharacterisation} with the function $x\mapsto\psi_j(\pm x)$ gives the asserted equivalence stated in the proposition.
\end{proof}

\section{Method of Differential Inequalities}
\label{sec:Method}

Throughout this section, we assume that $U$ is a unitary operator at least of class $C_{\mathrm u}^1(A)$. We denote the spectral measure of $U$ by $E$. Let $\zeta_0\in\mu^A(U)$. By the definition of the Mourre set, there exist an open subset $\mathbb{T} \neq I\subset\mu^A(U)$, containing $\zeta_0$, and a constant $a_I>0$ such that
\begin{subequations}\label{eq:MourreEstimateInterval}
	\begin{align}
		E(I)[A,U]U^*E(I)\geq a_I E(I).
		\label{eq:MourreEstimateInterval1}
	\end{align}
	Because $U$ commutes with $E(I)$, conjugation by $U$ also gives
	\begin{align}
		E(I)U^*[A,U]E(I)=U^*(E(I)[A,U]U^*E(I))U\geq a_I E(I).
		\label{eq:MourreEstimateInterval2}
	\end{align}
\end{subequations}
We fix a compact subset $I_0\subset I$. Note that $d=\operatorname{dist}(I_0,\mathbb T\backslash I)>0$ because $I$ is open and $I_0$ is a compact subset.

In the following, we write a complex spectral parameter as $z=r\zeta$, where $r\in[0,\infty)$ and $\zeta\in I_0$. To treat the interior and exterior of the unit circle simultaneously, we choose $\sigma\in\{-1,1\}$ so that $\sigma(1-r^2)=|1-r^2|$. When $r=1$, either choice for $\sigma$ is allowed.

\subsection{Basic Estimates}
\label{ssec:BasicEstimates}

The proof of the following lemma is the only point at which we apply the Mourre estimate~\eqref{eq:MourreEstimateInterval} directly.

\begin{lem}\label{lem:BasicEstimate}
	For every $a\in(0,a_I)$, there exist $\varepsilon_0>0$ and a constant $c<\infty$ such that, for all $r\in[0,\infty)$, $\zeta\in I_0$, $\varepsilon\in[0,\varepsilon_0]$, and $f\in\mathcal H$,
	\begin{subequations}\label{eq:BasicEstimate}
		\begin{align}
			(2ar^2\varepsilon+|1-r^2|)\|f\|^2 &\leq |\scp{f}{(1-r^2U_{\sigma\varepsilon}(U_{\sigma\varepsilon})^*)f}|
			+ c\varepsilon \|(1-r\zeta (U_{\sigma\varepsilon})^*)f\|^2, \label{eq:BasicEstimate1} \\
			(2ar^2\varepsilon+|1-r^2|)\|f\|^2
			&\leq	|\scp{f}{(1-r^2(U_{\sigma\varepsilon})^*U_{\sigma\varepsilon})f}|
			+c\varepsilon\|(1-r\overline\zeta U_{\sigma\varepsilon})f\|^2. \label{eq:BasicEstimate2}
		\end{align}
	\end{subequations}
\end{lem}

\begin{proof}
	We prove only \eqref{eq:BasicEstimate1} because the proof of \eqref{eq:BasicEstimate2} is very similar. We set $P=E(I)$ and $P^\perp = 1-P$. To prove \eqref{eq:BasicEstimate1}, we apply the Mourre estimate on $P\mathcal{H}$ to control the first-order unitarity defect of $U_{\sigma\varepsilon}$, and, on the complementary spectral subspace $P^\perp \mathcal{H}$, we use the positive distance between $I_0$ and $\mathbb T\backslash I$.
	
	Let $f\in\mathcal{H}$. We first establish an estimate on $\|P^\perp f\|$. For $w\in\mathbb T\backslash I$, $r\in[0,\infty)$, and $\zeta\in I_0$, the triangle inequality gives
	\begin{align}
		rd\leq r|w-\zeta|\leq |rw-\zeta|+|1-r|\leq 2|rw-\zeta|=2|1-r\zeta\overline w|,
	\end{align}
	where $d>0$ is the distance between $I_0$ and $\mathbb{T}\backslash I$. This implies that
	\begin{align}
		r^2\|P^\perp f\|^2 \leq \frac{4}{d^2}\|(1-r\zeta U^*)f\|^2.
	\end{align}
	To replace $U$ with $U_{\sigma\varepsilon}$, we use that, by Lemma~\ref{lem:Uepsilon}, $\|(U_{\sigma\varepsilon})^*-U^*\|\leq c_0 \varepsilon$. Thus, there exists a constant $c_1$, independent of $r\in[0,\infty)$, $\zeta\in I_0$, and $\varepsilon$, such that
	\begin{align} \label{eq:EstimatePperp}
		r^2\|P^\perp f\|^2 \leq c_1 \|(1-r\zeta (U_{\sigma\varepsilon})^*)f\|^2 + c_1 r^2\varepsilon^2\|f\|^2.		
	\end{align}	
	Next, we estimate the unitary defect of $U_{\sigma\varepsilon}$. Choose $a_0\in (a,a_I)$. Note that the Mourre estimate \eqref{eq:MourreEstimateInterval1} may be written $-P\mathcal{A}[U]U^*P \geq a_I P$. Using \eqref{eq:ConvergenceToCommutator1}, we see that there exists $\varepsilon_1>0$ such that, for $\varepsilon\in [0,\varepsilon_1]$,
	\begin{align}
		\|\sigma(1 - U_{\sigma \varepsilon} (U_{\sigma \varepsilon})^*) + 2 \varepsilon \mathcal{A}[U]U^*\| \leq 2\varepsilon(a_I - a_0).
	\end{align}
	Then, from the Mourre estimate~\eqref{eq:MourreEstimateInterval1}, it follows that, for $\varepsilon\in[0,\varepsilon_1]$,
	\begin{align} \label{eq:MourreEstimateIntervalEpsilon}
		2 a_0 \varepsilon P \leq \sigma P(1 - U_{\sigma \varepsilon} (U_{\sigma \varepsilon})^*)P.
	\end{align}		
	Moreover, by \eqref{eq:ConvergenceToCommutator1}, it holds that $\|1 - U_{\sigma \varepsilon} (U_{\sigma \varepsilon})^*\| \leq c_2 \varepsilon$. Using \eqref{eq:MourreEstimateIntervalEpsilon}, we obtain
	\begin{align}
		&\sigma \scp{f}{(1 - U_{\sigma \varepsilon} (U_{\sigma \varepsilon})^*)f}
		= \sigma \scp{Pf}{(1 - U_{\sigma \varepsilon} (U_{\sigma \varepsilon})^*)Pf} \nonumber \\
		&\qquad +2\sigma \Re \scp{P^\perp f}{(1 - U_{\sigma \varepsilon} (U_{\sigma \varepsilon})^*)Pf} + \sigma \scp{P^\perp f}{(1 - U_{\sigma \varepsilon} (U_{\sigma \varepsilon})^*) P^\perp f} \nonumber \\
		&\quad \geq 2a_0\varepsilon\|Pf\|^2 - 2c_2\varepsilon\|Pf\|\|P^\perp f\| - c_2\varepsilon\|P^\perp f\|^2.
	\end{align}
	To estimate the mixed term, we apply Young's inequality:
	\begin{align}
		2c_2\varepsilon\|Pf\|\|P^\perp f\| \leq (a_0-a)\varepsilon\|Pf\|^2+  \frac{c_2^2\varepsilon}{a_0-a}\|P^\perp f\|^2,
	\end{align}
	where $a_0-a>0$. Thus, using $\|Pf\|^2 = \|f\|^2 - \|P^\perp f\|^2$, we see that there exists a constant $c_3$ such that
	\begin{align}
		\sigma \scp{f}{(1 - U_{\sigma \varepsilon} (U_{\sigma \varepsilon})^*)f} \geq
		(a+a_0)\varepsilon\|f\|^2 - c_3 \varepsilon \|P^\perp f\|^2.
	\end{align}	
	Combining the preceding estimate with \eqref{eq:EstimatePperp}, we find that
	\begin{align}
		&\sigma \scp{f}{(1 - r^2 U_{\sigma \varepsilon} (U_{\sigma \varepsilon})^*) f} 
		= |1-r^2| \|f\|^2 + \sigma r^2 \scp{f}{(1 - U_{\sigma \varepsilon} (U_{\sigma \varepsilon})^*) f} 
		\nonumber \\
		&\quad \geq (|1-r^2|+(a+a_0)r^2\varepsilon-c_1c_3r^2\varepsilon^3)\|f\|^2 - c_1 c_3 \varepsilon\|(1-r\zeta (U_{\sigma\varepsilon})^*)f\|^2.
	\end{align}
	Finally, we choose $\varepsilon_0\in(0,\varepsilon_1]$ such that $c_1c_3\varepsilon_0^2 \leq a_0-a$. Then, for $\varepsilon\in[0,\varepsilon_0]$, we obtain \eqref{eq:BasicEstimate1}.
\end{proof}

\subsection{Regularised Resolvents}

From the estimates of Lemma~\ref{lem:BasicEstimate}, we deduce that $1-r\zeta(U_{\sigma\varepsilon})^*$ is invertible for $\varepsilon\in (0,\varepsilon_0]$, where $r$ and $\zeta$ are as above. For $r\ne1$, this is obvious because $(U_{\sigma\varepsilon})^*$ converges to $U^*$ in norm as $\varepsilon\to0$. The non-trivial part is that $1-r\zeta(U_{\sigma\varepsilon})^*$ is invertible also in the case $r=1$ and $\varepsilon\in (0,\varepsilon_0]$. 

\begin{prop}\label{prop:QuadraticEstimate}
	Let $a$, $\varepsilon_0$, and $c$ be as in Lemma~\ref{lem:BasicEstimate}. Let $\varepsilon\in[0,\varepsilon_0]$, $r\in [0,\infty)$, and $\zeta\in I_0$. Assume that $r\neq 1$ or $\varepsilon>0$. Then, the operator $1-r\zeta (U_{\sigma\varepsilon})^*$ is invertible in $B(\mathcal H)$. Set
	\begin{align}
		G_{\sigma\varepsilon} &= G_{\sigma\varepsilon}(r,\zeta) = (1-r\zeta (U_{\sigma\varepsilon})^*)^{-1},\\
		H_{\sigma\varepsilon} &= H_{\sigma\varepsilon}(r,\zeta) = (1-r\zeta (U_{\sigma\varepsilon})^*)^{-1}(1+r\zeta (U_{\sigma\varepsilon})^*) = 2G_{\sigma\varepsilon}(r,\zeta)-1.
	\end{align}
	For every $f\in\mathcal H$, we have
	\begin{align} \label{eq:QuadraticEstimate}
		\|G_{\sigma\varepsilon}^{(*)}f\|^2 \leq \frac{|\Re\scp{f}{H_{\sigma\varepsilon}f}|+ c\varepsilon\|f\|^2}{2ar^2\varepsilon+|1-r^2|},		
	\end{align}
	where $G_{\sigma\varepsilon}^{(*)}$ denotes either $G_{\sigma\varepsilon}$ or $G_{\sigma\varepsilon}^*$. Moreover,
	\begin{align} \label{eq:QuadraticEstimate2}
		\|G_{\sigma\varepsilon}^{(*)}\| \leq \frac{1+r\|U_{\sigma\varepsilon}\|}{2ar^2\varepsilon+|1-r^2|}	+
		\left(\frac{c\varepsilon}{2ar^2\varepsilon+|1-r^2|}\right)^\frac{1}{2}.
	\end{align}
\end{prop}

\begin{proof}
	We first establish invertibility of the operator $1-r\zeta(U_{\sigma\varepsilon})^*$. The identities
	\begin{subequations}\label{eq:Identity}
		\begin{align}
			1-r^2U_{\sigma\varepsilon}(U_{\sigma\varepsilon})^*
			&=(1-r\zeta (U_{\sigma\varepsilon})^*)^*+r\overline\zeta U_{\sigma\varepsilon}(1-r\zeta (U_{\sigma\varepsilon})^*),	\label{eq:Identity1}\\
			1-r^2(U_{\sigma\varepsilon})^*U_{\sigma\varepsilon}
			&=(1-r\overline\zeta U_{\sigma\varepsilon})^*+r\zeta (U_{\sigma\varepsilon})^*(1-r\overline\zeta U_{\sigma\varepsilon})	\label{eq:Identity2},
		\end{align}
	\end{subequations}
	and Lemma~\ref{lem:BasicEstimate} imply that, for every $f\in\mathcal H$,
	\begin{subequations}
		\begin{align}
			(2ar^2\varepsilon+|1-r^2|)\|f\|	&\leq (1+r\|U_{\sigma\varepsilon}\|+c\varepsilon\|1-r\zeta (U_{\sigma\varepsilon})^*\|) \|(1-r\zeta (U_{\sigma\varepsilon})^*)f\|, \\
			(2ar^2\varepsilon+|1-r^2|)\|f\|	&\leq (1+r\|U_{\sigma\varepsilon}\|+c\varepsilon\|1-r\overline\zeta U_{\sigma\varepsilon}\|) \|(1-r\overline\zeta U_{\sigma\varepsilon})f\|.
		\end{align}
	\end{subequations}	
	Because $r\neq 1$ or $\varepsilon>0$, the coefficient $2ar^2\varepsilon+|1-r^2|$ is strictly positive. Thus, $1-r\zeta (U_{\sigma\varepsilon})^*$ and its adjoint $1-r\overline{\zeta}U_{\sigma\varepsilon}$ are bounded from below by a positive constant. In particular, $1-r\zeta (U_{\sigma\varepsilon})^*$ is injective and has closed range. Its range is also dense, because its adjoint has trivial kernel. Hence, the operator is bijective and has a bounded inverse.
	
	To prove \eqref{eq:QuadraticEstimate}, we replace $f$ in \eqref{eq:BasicEstimate} with $G_{\sigma\varepsilon}f$ and $G_{\sigma\varepsilon}^*f$, respectively. Indeed, using \eqref{eq:Identity}, we obtain
	\begin{subequations}
		\begin{align}
			G_{\sigma\varepsilon}^*	(1-r^2U_{\sigma\varepsilon}(U_{\sigma\varepsilon})^*)G_{\sigma\varepsilon}
			&=G_{\sigma\varepsilon}+G_{\sigma\varepsilon}^*-1 = 2\Re G_{\sigma\varepsilon} - 1,\\
			G_{\sigma\varepsilon}(1-r^2 (U_{\sigma\varepsilon})^*U_{\sigma\varepsilon})G_{\sigma\varepsilon}^*
			&=G_{\sigma\varepsilon}+G_{\sigma\varepsilon}^*-1 = 2\Re G_{\sigma\varepsilon} - 1.
		\end{align}
	\end{subequations}	
	Finally, because
	\begin{align}
		\|H_{\sigma\varepsilon}\| \leq (1+r\|U_{\sigma\varepsilon}\|) \|G_{\sigma\varepsilon}^{(*)}\|,
	\end{align}
	the estimate \eqref{eq:QuadraticEstimate} implies that
	\begin{align}
		\|G_{\sigma\varepsilon}^{(*)}\|^2 \leq \frac{(1+r\|U_{\sigma\varepsilon}\|)\|G_{\sigma\varepsilon}^{(*)}\| + c\varepsilon}{2ar^2\varepsilon+|1-r^2|},
	\end{align}
	which has \eqref{eq:QuadraticEstimate2} as a consequence.
\end{proof}

For fixed $\varepsilon\in(0,\varepsilon_0]$, Proposition~\ref{prop:QuadraticEstimate} gives that $1-r\zeta (U_{\sigma\varepsilon})^*$ is invertible at $r=1$ for both choices of $\sigma \in\{-1,1\}$. The two different choices of $\sigma$ give continuous extensions of the function $r \mapsto G_{\sigma\varepsilon}(r,\zeta)$, as $r$ approaches 1 from below or above, respectively.

\begin{cor}\label{cor:ExistenceLimitEpsilonPositive}
	Let $\varepsilon\in(0,\varepsilon_0]$. The functions
	\begin{subequations}
		\begin{align}
			[0,1)\times I_0 \ni (r,\zeta) & \mapsto G_{\varepsilon}(r,\zeta) \in B(\mathcal{H}), \\
			(1,\infty)\times I_0 \ni (r,\zeta) & \mapsto G_{-\varepsilon}(r,\zeta) \in B(\mathcal{H})
		\end{align}
	\end{subequations}
	have unique norm-continuous extensions to $[0,1]\times I_0$ and $[1,\infty)\times I_0$, respectively.
\end{cor}

\begin{proof}
	By Proposition~\ref{prop:QuadraticEstimate}, $G_{\varepsilon}(r,\zeta)$ is defined for $(r,\zeta)\in[0,1]\times I_0$, and $G_{-\varepsilon}(r,\zeta)$ is defined for $(r,\zeta)\in[1,\infty)\times I_0$. For either sign, the resolvent identity gives
	\begin{align}
		G_{\pm\varepsilon}(r_1,\zeta_1)-G_{\pm\varepsilon}(r_2,\zeta_2) 
		=(r_1\zeta_1-r_2\zeta_2)G_{\pm\varepsilon}(r_1,\zeta_1)	(U_{\pm\varepsilon})^*G_{\pm\varepsilon}(r_2,\zeta_2),
	\end{align}
	where $r_1$ and $r_2$ are assumed to be both in $[0,1]$ in the $+$-case, and in $[1,\infty)$ in the $-$-case. This implies the claim if we use the bounds on the operator norm of $G_{\pm\varepsilon}$ from \eqref{eq:QuadraticEstimate2}.
\end{proof}

\subsection{Differential Inequality}
\label{ssec:DifferentialInequality}

Let $\varphi\in C_0^\infty(\mathbb{R})$ be as in Section~\ref{sec:Regularisation}, that is, $\varphi$ is real-valued and $\varphi(x)=1$ in a neighbourhood of 0. For $f\in\mathcal{H}$, we define $f_\varepsilon = \varphi(\varepsilon A) f$. Then, $f_\varepsilon\to f$ in $\mathcal H$ as $\varepsilon\to0$, the map $\varepsilon\mapsto f_\varepsilon$ is of class $C^\infty$ on $\mathbb R\backslash\{0\}$, and $f_\varepsilon\in D(A^\infty)$ for every $\varepsilon\ne0$.

Let $r\in [0,\infty)$, $\zeta\in I_0$, and $\varepsilon\in [0,\varepsilon_0]$, where $\varepsilon_0$ is as in Lemma~\ref{lem:BasicEstimate}. Assume that $r\neq 1$ or $\varepsilon\neq 0$. According to Proposition~\ref{prop:QuadraticEstimate}, we can define the function
\begin{align}
	F_{\sigma\varepsilon} = F_{\sigma\varepsilon}(r,\zeta) = \scp{f_\varepsilon}{H_{\sigma\varepsilon}(r,\zeta) f_\varepsilon}
	= 2 \scp{f_\varepsilon}{G_{\sigma\varepsilon}(r,\zeta) f_\varepsilon} - \|f_\varepsilon\|^2.
\end{align}
The following identity is the starting point for the differential inequality that we derive in the proof of Proposition~\ref{prop:MajorisingFunction}. We denote derivatives with respect to $\varepsilon$ by a prime.

\begin{lem} \label{lem:RegularityFepsilon}
	For every $r\in[0,\infty)$ and $\zeta\in I_0$, the function $\varepsilon\mapsto F_{\sigma\varepsilon}(r,\zeta)$ is of class $C^\infty$ on $(0,\varepsilon_0]$, and, if $r\neq 1$, continuous on $[0,\varepsilon_0]$. Moreover, for every $\varepsilon\in (0,\varepsilon_0)$, its derivative is given by
	\begin{align} \label{eq:FepsilonPrime}
		\frac{1}{2} F_{\sigma\varepsilon}' = \scp{f_\varepsilon'+\sigma Af_\varepsilon}{G_{\sigma\varepsilon} f_\varepsilon} + \scp{G_{\sigma\varepsilon}^* f_\varepsilon}{f_\varepsilon' -\sigma Af_\varepsilon} 
		- \sigma r\zeta \scp{K_{\sigma\varepsilon} G_{\sigma\varepsilon}^* f_\varepsilon}{G_{\sigma\varepsilon} f_\varepsilon} - \Re \scp{f_\varepsilon'}{f_\varepsilon}.
	\end{align}
\end{lem}

\begin{proof}
	It is not difficult to verify that $\varepsilon\mapsto F_{\sigma\varepsilon}(r,\zeta)$ is of class $C^\infty$ on $(0,\varepsilon_0)$ because $U_{\sigma\varepsilon}$ and~$f_\varepsilon$ are of class $C^\infty$ on $(0,\varepsilon_0)$. It is also easy to see that, if $r\neq 1$, then $F_{\sigma\varepsilon}$ converges to $F_0 = 2\scp{f}{(1-r\zeta U^*)^{-1}f} - \|f\|^2$ as $\varepsilon\downarrow0$. Moreover, the derivative of $G_{\sigma\varepsilon} = (1-r\zeta (U_{\sigma\varepsilon})^*)^{-1}$ with respect to $\varepsilon$ is given by
	\begin{align} \label{eq:DerivativeG}
		\partial_\varepsilon G_{\sigma\varepsilon} &= G_{\sigma\varepsilon} r \zeta \partial_\varepsilon (U_{\sigma\varepsilon})^* G_{\sigma\varepsilon} \nonumber \\
		&= -\sigma G_{\sigma\varepsilon} r\zeta (\mathcal{A}[(U_{\sigma\varepsilon})^*] + K_{\sigma\varepsilon}^*) G_{\sigma\varepsilon} = -\sigma \mathcal{A}[G_{\sigma\varepsilon}] - \sigma r\zeta G_{\sigma\varepsilon} K_{\sigma\varepsilon}^* G_{\sigma\varepsilon},
	\end{align}
	where we used \eqref{eq:UepsilonPrime} and that $G_{\sigma\varepsilon}\in C^1(A)$ \cite[Proposition~5.1.6]{amrein1996}. This gives the claimed formula for $F_{\sigma\varepsilon}'$. 
\end{proof}

For $\varepsilon\in(0,\varepsilon_0]$, the functions $(r,\zeta) \mapsto F_{\sigma\varepsilon}(r,\zeta)$ and their derivatives with respect to $\varepsilon$ remain continuous as $r$ approaches 1 from either below or above. This follows from Corollary~\ref{cor:ExistenceLimitEpsilonPositive} and the derivative formula in Lemma~\ref{lem:RegularityFepsilon}.

\begin{cor} \label{cor:ExtensionsFepsilon}	
	Let $\varepsilon\in(0,\varepsilon_0]$. For $k\in\{0,1\}$, the functions
	\begin{subequations}
		\begin{align}
			[0,1)\times I_0 \ni (r,\zeta) & \mapsto F_\varepsilon^{(k)}(r,\zeta), \\
			(1,\infty)\times I_0 \ni (r,\zeta) & \mapsto F_{-\varepsilon}^{(k)}(r,\zeta)
		\end{align}
	\end{subequations}
	have unique continuous extensions to $[0,1]\times I_0$ and $[1,\infty)\times I_0$, respectively.
\end{cor}

The next proposition provides an integrable bound for the function $\varepsilon\mapsto F_{\sigma\varepsilon}'(r,\zeta)$, uniformly in $r\in [0,\infty)$ and $\zeta\in I_0$. This bound is the main result of the method of differential inequalities. To obtain this bound, we need to use the regularity assumptions on $U$ and $f$ imposed by Theorem~\ref{thm:LAP}.

\begin{prop}\label{prop:MajorisingFunction}
	Assume that $U\in\mathcal C^{1,1}(A)$ and $f\in\mathcal{H}_{1/2,1}$. Then, there exists a function $\rho\in L^1((0,\varepsilon_0))$ such that, for every $\varepsilon\in (0,\varepsilon_0)$,
	\begin{align}
		\sup_{r\in[0,\infty),\,\zeta\in I_0} |F_{\sigma\varepsilon}'(r,\zeta)| \leq \rho(\varepsilon).
	\end{align}
\end{prop}

\begin{proof}
	Without loss of generality, we may assume that $2a\varepsilon_0 \leq 1$ because only small $\varepsilon$ are critical. Indeed, in the region where $\varepsilon$ is strictly separated from 0, we can use \eqref{eq:QuadraticEstimate2} to bound $F_{\sigma\varepsilon}'(r,\zeta)$.
	
	We set $\ell(\varepsilon) = \|f_\varepsilon'\| + \|Af_\varepsilon\|$. Then, from \eqref{eq:FepsilonPrime}, it follows that
	\begin{align} \label{eq:DifferentialInequality}
		\frac{1}{2} |F_{\sigma\varepsilon}'| \leq \ell(\varepsilon) (\|G_{\sigma\varepsilon}f_\varepsilon\|
		+ \|G_{\sigma\varepsilon}^*f_\varepsilon\|) + r\|K_{\sigma\varepsilon}\| \|G_{\sigma\varepsilon}f_\varepsilon\|
		\|G_{\sigma\varepsilon}^*f_\varepsilon\| + \|f_\varepsilon'\|\|f_\varepsilon\|.
	\end{align}
	We bound the factors $\|G_{\sigma\varepsilon}^{(*)}f_\varepsilon\|$ and $r\|G_{\sigma\varepsilon}^{(*)}f_\varepsilon\|$ by \eqref{eq:QuadraticEstimate}. To simplify the estimate, we use that $2ar^2\varepsilon+|1-r^2| \geq 2a\varepsilon \max\{1,r^2\}$ for $\varepsilon\in [0,\varepsilon_0]$. This is readily verified by considering the cases $r\in [0,1)$, $r=1$, and $r\in (1,\infty)$, separately. It makes use of the assumption $2a\varepsilon_0\leq 1$. We obtain the following estimates:
	\begin{subequations} \label{eq:EstiamteG}
		\begin{align}
			\|G_{\sigma\varepsilon}^{(*)}f_\varepsilon\| 
			&\leq \frac{1}{\sqrt{2a\varepsilon}} |F_{\sigma\varepsilon}|^\frac{1}{2} + \sqrt{\frac{c}{2a}} \|f_\varepsilon\|, \\
			r \|G_{\sigma\varepsilon}^{(*)}f_\varepsilon\| 
			&\leq \frac{1}{\sqrt{2a\varepsilon}} |F_{\sigma\varepsilon}|^\frac{1}{2} + \sqrt{\frac{c}{2a}} \|f_\varepsilon\|.
		\end{align}
	\end{subequations}
	We insert the estimates \eqref{eq:EstiamteG} into \eqref{eq:DifferentialInequality}. It follows that there exists a constant $b>0$, independent of $r\in [0,\infty)$, $\zeta\in I_0$, and $\varepsilon \in (0,\varepsilon_0)$, such that
	\begin{align} \label{eq:BoundFepsilonPrime}
		b|F_{\sigma\varepsilon}'| \leq \ell(\varepsilon)\|f_\varepsilon\| + \|K_{\sigma\varepsilon}\| \|f_\varepsilon\|^2
		+\frac{\ell(\varepsilon)}{\sqrt{\varepsilon}} |F_{\sigma\varepsilon}|^\frac{1}{2}	+\frac{\|K_{\sigma\varepsilon}\|}{\varepsilon}|F_{\sigma\varepsilon}|.
	\end{align}	
	This is the differential inequality, which bounds the derivative of $F_{\sigma\varepsilon}$ by the function itself. From Grönwall's inequality \cite[Appendix~7.A, (7.3.11)]{amrein1996} \cite[p.~41]{boutetDeMonvel1996}, we obtain, for
	$\varepsilon\in(0,\varepsilon_0)$, the following bound on $F_{\sigma\varepsilon}$:
	\begin{align}
		|F_{\sigma\varepsilon}|	\leq 2\left\{ |F_{\sigma\varepsilon_0}| + \frac{1}{b}\int_\varepsilon^{\varepsilon_0} [\ell(\tau)\|f_\tau\| + \|K_{\sigma\tau}\| \|f_\tau\|^2] \diff \tau
		+ \frac{1}{b^2} \left(\int_\varepsilon^{\varepsilon_0}\frac{\ell(\tau)}{\sqrt{\tau}} \diff\tau \right)^2 \right\} \nonumber \\
		\exp\left( \frac{1}{b} \int_\varepsilon^{\varepsilon_0}\|K_{\sigma\tau}\|\frac{\diff\tau}{\tau} \right).
	\end{align}	
	We analyse the convergence of the integrals as $\varepsilon\downarrow 0$. By Proposition~\ref{prop:Regularity}, $U\in \mathcal{C}^{1,1}(A)$ if and only if
	\begin{align}
		\int_0^{\varepsilon_0} \|K_{\sigma\tau}\| \frac{\diff\tau}{\tau}<\infty.
	\end{align}
	Furthermore, let $f=f_0+f_1$ be a decomposition of $f$ with $f_0\in\mathcal{H}$ and $f_1\in D(A)$. Then,
	\begin{align}
		\ell(\varepsilon) \leq \varepsilon^{-1} (\|x\varphi'\|_{L^\infty} + \|x\varphi\|_{L^\infty}) \|f_0\| + (\|\varphi'\|_{L^\infty} + \|\varphi\|_{L^\infty}) \|Af_1\|.
	\end{align}
	Taking the infimum over all such decompositions of $f$, we see that there exists a constant $C<\infty$, depending only on $\varphi$, such that $\ell(\varepsilon) \leq C \varepsilon^{-1} K(\varepsilon,f)$, where $K(\varepsilon,\cdot)$ is the $K$-functional for the couple $(D(A),\mathcal{H})$, as defined in \cite[(2.2.1)]{amrein1996}. From the definition of the interpolation space $\mathcal{H}_{1/2,1} = (D(A),\mathcal{H})_{1/2,1}$ \cite[Definition~2.2.3]{amrein1996}, it follows that
	\begin{align}
		\int_0^{\varepsilon_0}\frac{ \ell(\tau)}{\sqrt{\tau}} \diff\tau 
		\leq C \int_0^{\varepsilon_0} \tau^{-\frac{1}{2}} K(\tau,f) \frac{\diff \tau}{\tau} \leq C \|f\|_{1/2,1}.
	\end{align}	
	Moreover, $F_{\sigma\varepsilon_0}(r,\zeta)$ is uniformly bounded in $r\in [0,\infty)$ and $\zeta\in I_0$ by Proposition~\ref{prop:QuadraticEstimate}. To summarise the above steps, we have proved that
	\begin{align} \label{eq:BoundOnF}
		\sup_{\varepsilon\in (0,\varepsilon_0)} \sup_{r\in[0,\infty),\,\zeta\in I_0} |F_{\sigma\varepsilon}(r,\zeta)| \leq C \|f\|^2_{1/2,1}.
	\end{align}
	Substituting this bound into \eqref{eq:BoundFepsilonPrime} gives the claim of the proposition.
\end{proof}

\begin{cor} \label{cor:H0sEstimate}
	Assume that $U\in\mathcal C^{1,1}(A)$. Then, there exists a constant $C<\infty$ such that, for all $\varepsilon\in (0,\varepsilon_0)$,
	\begin{align} \label{eq:EstimateGH121}
		\sup_{r\in[0,\infty),\,\zeta\in I_0} \left( \|G_{\sigma\varepsilon}(r,\zeta)\|_{\mathcal{B}(\mathcal{H}_{1/2,1},\mathcal{H})} + \|G_{\sigma\varepsilon}(r,\zeta)\|_{\mathcal{B}(\mathcal{H},\mathcal{H}_{-1/2,\infty})}\right) \leq C \varepsilon^{-\frac{1}{2}}.
	\end{align}	
	In particular, for every $s>1/2$, there exists a constant $C<\infty$ such that
	\begin{align}
		\sup_{r\in[0,\infty),\,\zeta\in I_0} \left( \|G_{\sigma\varepsilon}(r,\zeta)\langle A\rangle^{-s}\| +\|\langle A\rangle^{-s}G_{\sigma\varepsilon}(r,\zeta)\| \right) \leq C \varepsilon^{-\frac{1}{2}}.
	\end{align}
\end{cor}

\begin{proof}
	As before, we may assume that $2a\varepsilon_0\leq 1$. Because $\mathcal{H}_{-1/2,\infty}$ is the dual of $\mathcal{H}_{1/2,1}$, the polarisation method	applied to \eqref{eq:BoundOnF} gives
	\begin{align} \label{eq:EstimateGH121a}
		\sup_{\varepsilon\in (0,\varepsilon_0)} \sup_{r\in[0,\infty),\,\zeta\in I_0} \|\varphi(\varepsilon A) G_{\sigma\varepsilon}(r,\zeta) \varphi(\varepsilon A) \|_{\mathcal{B}(\mathcal{H}_{1/2,1},\mathcal{H}_{-1/2,\infty})} < \infty.
	\end{align}
	Let $f\in\mathcal{H}_{1/2,1}$. Then, applying \eqref{eq:EstiamteG} to $\varphi(\varepsilon A)f$, we obtain that
	\begin{align} \label{eq:EstimateG2}
		\|G_{\sigma\varepsilon}^{(*)} \varphi(\varepsilon A) f\| &\leq \frac{1}{\sqrt{2a\varepsilon}} |\scp{f}{\varphi(\varepsilon A)(2G_{\sigma\varepsilon}-1)\varphi(\varepsilon A)f}|^{\frac{1}{2}} + \sqrt{\frac{c}{2a}} \|\varphi\|_{L^\infty} \|f\| \nonumber \\
		&\leq C \varepsilon^{-\frac{1}{2}} \|f\|_{\mathcal{H}_{1/2,1}},
	\end{align}
	where the second inequality follows from \eqref{eq:EstimateGH121a}
	and the continuous embedding $\mathcal{H}_{1/2,1}\subset\mathcal{H}$. To remove $\varphi(\varepsilon A)$ in the estimate, note that $|1-\varphi(x)|\leq C|x|^{1/2}$ because $\varphi=1$ near zero and $\varphi$ is bounded. This implies that
	\begin{align}
		\|G_{\sigma\varepsilon}^{(*)}(1-\varphi(\varepsilon A))f\|
		\leq C\varepsilon^{-1}\|(1-\varphi(\varepsilon A))f\|
		\leq C\varepsilon^{-\frac12}\|f\|_{\mathcal{H}_{1/2}}
		\leq C\varepsilon^{-\frac12}\|f\|_{\mathcal{H}_{1/2,1}},
	\end{align}
	where we used the bound $\|G_{\sigma\varepsilon}^{(*)}\|\leq C\varepsilon^{-1}$ from \eqref{eq:QuadraticEstimate2} and the continuous embedding $\mathcal{H}_{1/2,1}\subset\mathcal{H}_{1/2}$. Combining the last estimate with \eqref{eq:EstimateG2} proves \eqref{eq:EstimateGH121}. 
	
	The final assertion of the corollary follows from $\mathcal{H}_s\subset\mathcal{H}_{1/2,1}$ for every $s>1/2$.
\end{proof}

\subsection{Higher-Order Estimates}
\label{ssec:HigherOrderEstimates}

In this subsection, we establish estimates on the derivatives of the regularised resolvent, which are needed in the proof of Theorem~\ref{thm:HigherRegularity}.

Let $\varphi\in C_0^\infty(\mathbb{R})$ be a real-valued function as above, that is, $\varphi=1$ in a neighbourhood of 0. For $\zeta\in I_0$ and $\varepsilon\in (0,\varepsilon_0)$, we define the functions
\begin{align}
	L_{\pm}(\zeta,\varepsilon) = \varphi(\varepsilon A) G_{\pm\varepsilon}(1,\zeta) \varphi(\varepsilon A).
\end{align}
The goal of this subsection is to control the derivatives in $\zeta$ and $\varepsilon$ of this function. We recall the following two estimates from Proposition~\ref{prop:QuadraticEstimate} and Corollary~\ref{cor:H0sEstimate}, which hold for $\varepsilon\in (0,\varepsilon_0)$:
\begin{align}
	&\sup_{\zeta\in I_0} \|G_{\pm\varepsilon}(1,\zeta)\| \leq C \varepsilon^{-1}, \label{eq:SimpleEstimateG} \\
	&\sup_{\zeta\in I_0} \left( \|G_{\pm\varepsilon}(1,\zeta)\langle A\rangle^{-s}\| +\|\langle A\rangle^{-s}G_{\pm\varepsilon}(1,\zeta)\| \right) \leq C \varepsilon^{-\frac{1}{2}} \label{eq:EstimateAG}.
\end{align}
The first estimate requires $U\in C^1_\mathrm{u}(A)$, the second $U\in \mathcal{C}^{1,1}(A)$. If $U\in \mathcal{C}^{s+1/2}(A)$ with $s>1/2$ as in Theorem~\ref{thm:HigherRegularity}, then both assumptions are satisfied because $\mathcal{C}^{s+1/2}(A) \subset \mathcal{C}^{1,1}(A) \subset C_\mathrm{u}^1(A)$. 

We start with the following lemma estimating the derivatives of $G_{\pm\varepsilon}(1,\zeta)$. This lemma will be useful to bound the derivatives of $L_{\pm}$.

\begin{lem} \label{lem:DerivativesG}
	Let $s>1/2$, and assume that $U\in \mathcal{C}^{1,1}(A)$. Then, for all integers $\ell,k\geq0$, there exists a constant $C<\infty$ such that, for all $\varepsilon\in (0,\varepsilon_0)$,
	\begin{align}\label{eq:GlobalOneSidedMixedDerivatives}
		\sup_{\zeta\in I_0} \left( \|D_\zeta^\ell\partial_\varepsilon^k G_{\sigma\varepsilon}(1,\zeta) \langle A\rangle^{-s}\|
		+\|\langle A\rangle^{-s} D_\zeta^\ell\partial_\varepsilon^k G_{\sigma\varepsilon}(1,\zeta)\| \right)
		\leq C \varepsilon^{-\ell-k-\frac{1}{2}}.
	\end{align}
\end{lem}

\begin{proof}
	In the first step, we estimate the derivatives of $1-\zeta (U_{\sigma\varepsilon})^*$. From Lemma~\ref{lem:Uepsilon}, it follows that, for all integers $p,q\geq 0$ with $p+q\geq1$, there exists a constant $C<\infty$ such that
	\begin{align}\label{eq:MixedLemmaInverseFactor}
		\| D_\zeta^p\partial_\varepsilon^q (1-\zeta (U_{\pm\varepsilon})^*) \| \leq C \varepsilon^{1-p-q}.
	\end{align}
	The power of $\varepsilon$ is not optimal in $p$ because derivatives in $\zeta$ do not increase the singularity at $\varepsilon=0$. However, stating the estimate in this way simplifies the power counting below.
	
	In the second step, we compute the derivatives of $G_{\pm\varepsilon}(1,\zeta)$. Differentiating $G_{\pm\varepsilon}(1,\zeta)$ once with respect to either $\zeta$ or $\varepsilon$ gives
	\begin{subequations}
		\begin{align}
			D_\zeta G_{\pm\varepsilon} &= -G_{\pm\varepsilon} (D_\zeta(1-\zeta (U_{\pm\varepsilon})^*)) G_{\pm\varepsilon}, \\
			\partial_\varepsilon G_{\pm\varepsilon} &= -G_{\pm\varepsilon} (\partial_\varepsilon(1-\zeta (U_{\pm\varepsilon})^*)) G_{\pm\varepsilon}.
		\end{align}
	\end{subequations}	
	The limits defining the derivatives exist in the operator norm topology. By induction, $D_\zeta^\ell\partial_\varepsilon^k G_{\pm\varepsilon}$ is a finite linear combination of products of the form
	\begin{align}\label{eq:MixedLemmaProduct}
		G_{\pm\varepsilon} \prod_{j=1}^{n} \left[ (D_\zeta^{p_j} \partial_\varepsilon^{q_j} (1-\zeta (U_{\pm\varepsilon})^*) ) G_{\pm\varepsilon} \right],
	\end{align}
	where $1\leq n\leq\ell+k$, $p_j,q_j\geq0$, $p_j+q_j\geq1$, $\sum_{j=1}^{n}p_j=\ell$, and $\sum_{j=1}^{n}q_j=k$. If we multiply \eqref{eq:MixedLemmaProduct} on the right with $\langle A\rangle^{-s}$, we obtain the following norm estimate:
	\begin{align}
		&\bigg\|G_{\pm\varepsilon} \prod_{j=1}^{n} \left[ (D_\zeta^{p_j} \partial_\varepsilon^{q_j}	(1-\zeta (U_{\pm\varepsilon})^*) ) G_{\pm\varepsilon} \right] \angles{A}^{-s}\bigg\|
		\nonumber \\
		&\quad \leq \|G_{\pm\varepsilon}\|^n \prod_{j=1}^{n} \| D_\zeta^{p_j} \partial_\varepsilon^{q_j}	(1-\zeta (U_{\pm\varepsilon})^*) \| \|G_{\pm\varepsilon} \angles{A}^{-s}\| \nonumber \\
		&\quad \leq C \varepsilon^{-n}
		\left(\prod_{j=1}^{n} \varepsilon^{1-p_j-q_j} \right) \varepsilon^{-\frac{1}{2}}
		\leq C \varepsilon^{-\ell-k-\frac{1}{2}},
	\end{align}
	where we used \eqref{eq:SimpleEstimateG}, \eqref{eq:EstimateAG}, and \eqref{eq:MixedLemmaInverseFactor} in passing from the second line to the third. If we multiply \eqref{eq:MixedLemmaProduct} on the left with $\langle A\rangle^{-s}$, we obtain the same estimate by a similar argument.
\end{proof}

\begin{prop}\label{prop:MixedDerivativeEstimates}
	Let $s> 1/2$, and assume that $U\in\mathcal{C}^{s+1/2}(A)$. Let $m\in\mathbb{N}$ be such that $m\geq2s$. Then, for every integer $\ell\geq0$, there exists a constant $C<\infty$ such that, for all $\varepsilon\in (0,\varepsilon_0)$,
	\begin{align}\label{eq:HighestMixedDerivativeEstimate}
		\sup_{\zeta\in I_0} \|\langle A\rangle^{-s} D_\zeta^\ell\partial_\varepsilon^m L_{\pm}(\zeta,\varepsilon) \langle A\rangle^{-s} \| \leq C \varepsilon^{s-\frac{1}{2}-m-\ell}.
	\end{align}
	In particular,
	\begin{align} \label{eq:MixedDerivativeEstimate}
		\sup_{\zeta\in I_0} \sum_{\substack{\ell,k\geq 0 \\ \ell+k=m}} \| \langle A\rangle^{-s} D_\zeta^\ell\partial_\varepsilon^k	L_{\pm}(\zeta,\varepsilon) \langle A\rangle^{-s} \| \leq C \varepsilon^{s-\frac{1}{2}-m}.
	\end{align}
\end{prop}

\begin{proof}	
	For an integer $i\geq0$, we define the function $\varphi_i\in C_0^\infty(\mathbb{R})$ by $\varphi_i(x)=x^i\varphi^{(i)}(x)$. Using the Leibniz rule, the derivatives of $L_{\pm}$ are given by
	\begin{align}\label{eq:CutoffLeibnizFormula}
		D_\zeta^\ell\partial_\varepsilon^m	L_{\pm}
		= \sum_{\substack{i,j,k\geq 0 \\i+j+k=m}} \frac{m!}{i!j!k!} \varepsilon^{k-m} \varphi_i(\varepsilon A)
		\, D_\zeta^\ell\partial_\varepsilon^k G_{\pm\varepsilon} \, \varphi_j(\varepsilon A).
	\end{align}
	To prove the proposition, we must show that, for $i+j+k=m$,
	\begin{align}\label{eq:GlobalCutoffMixedEstimate}
		\| \langle A\rangle^{-s} \varphi_i(\varepsilon A) \, D_\zeta^\ell\partial_\varepsilon^k
		G_{\pm\varepsilon} \, \varphi_j(\varepsilon A) \langle A\rangle^{-s} \| 
		\leq C \varepsilon^{s-\frac{1}{2}-\ell-k}.
	\end{align}		
	If $k=0$, then \eqref{eq:GlobalCutoffMixedEstimate} follows immediately from \eqref{eq:SimpleEstimateG} and \eqref{eq:EstimateAG}. Indeed, in this case, we have $i+j=m\geq 2s$; hence, $i\geq s$ or $j\geq s$. Let us assume that $i\geq s$. The proof for $j\geq s$ is similar. The derivative of $G_{\pm\varepsilon}$ with respect to $\zeta$ equals
	\begin{align}\label{eq:TangentialResolventDerivative}
		D_\zeta G_{\pm\varepsilon} = \I G_{\pm\varepsilon} \zeta (U_{\pm\varepsilon})^* G_{\pm\varepsilon} = \I(G_{\pm\varepsilon}^2-G_{\pm\varepsilon}).
	\end{align}
	Thus, $D_\zeta^\ell G_{\pm\varepsilon}$ is a linear combination of $G_{\pm\varepsilon}^p$, where $1\leq p\leq\ell+1$. For each such an integer $p$, it holds that
	\begin{align}
		\| \langle A\rangle^{-s} \varphi_i(\varepsilon A) G_{\pm\varepsilon}^p \varphi_j(\varepsilon A) \langle A\rangle^{-s} \| 
		&\leq \|\angles{A}^{-s} \varphi_i(\varepsilon A)\| \|G_{\pm\varepsilon}\|^{p-1} \|G_{\pm\varepsilon} \angles{A}^{-s} \| \|\varphi_j\|_{L^\infty} \nonumber \\
		&\leq C \varepsilon^s \varepsilon^{-p+1} \varepsilon^{-\frac12} 
		\leq C \varepsilon^{s-\frac{1}{2}-\ell},
	\end{align}
	where we used that $\|\angles{A}^{-s} \varphi_i(\varepsilon A)\|\leq C \varepsilon^s$, which follows from $i\geq s$.
	
	If $k\geq 1$, we need to compute the derivative of $G_{\pm\varepsilon}$ with respect to $\varepsilon$. By \eqref{eq:DerivativeG}, the first derivative with respect to $\varepsilon$ is given by
	\begin{align}\label{eq:RegularisationResolventDerivative}
		\partial_\varepsilon G_{\pm\varepsilon} = \mp\mathcal A[G_{\pm\varepsilon}]
		\mp \zeta G_{\pm\varepsilon} K_{\pm\varepsilon}^* G_{\pm\varepsilon}.
	\end{align}
	Iterating this identity, we obtain the formula
	\begin{align}\label{eq:IteratedRegularisationDerivative}
		D_\zeta^\ell \partial_\varepsilon^kG_{\pm\varepsilon} = (\mp1)^k\mathcal A^k[D_\zeta^\ell G_{\pm\varepsilon}]
		+ \sum_{n=0}^{k-1}(\mp1)^{k-n} \mathcal A^{k-1-n} [ D_\zeta^\ell \partial_\varepsilon^n (\zeta G_{\pm\varepsilon} K_{\pm\varepsilon}^*G_{\pm\varepsilon})	].
	\end{align}	
	Let us consider the first term on the right-hand side of \eqref{eq:IteratedRegularisationDerivative}. For an integer $a\geq 0$, we define the function $\varphi_{i,a} \in C_0^\infty(\mathbb{R})$ by $\varphi_{i,a}(x) = x^a \varphi_i(x) = x^{a+i} \varphi^{(i)}(x)$. Expanding the commutators as in \cite[p.~22]{boutetDeMonvel1996}, we obtain the identity
	\begin{align}
		\varphi_i(\varepsilon A)\mathcal A^k[D_\zeta^\ell G_{\pm\varepsilon}] \varphi_j(\varepsilon A) = \sum_{\substack{a,b\geq 0 \\ a+b=k}} \frac{k!}{a!b!} (-1)^a \varepsilon^{-k} \varphi_{i,a}(\varepsilon A) D_\zeta^\ell G_{\pm\varepsilon} \varphi_{j,b}(\varepsilon A).
	\end{align}
	We have $(a+i)+(b+j) = i + j + k = m \geq 2s$; hence, $a+i\geq s$ or $b+j\geq s$. Let us assume that $a+i\geq s$. The proof for $b+j\geq s$ is similar. Again, we write $D_\zeta^\ell G_{\pm\varepsilon}$ as a linear combination of $G_{\pm\varepsilon}^p$, where $1\leq p\leq\ell+1$. For each such an integer $p$, it holds that
	\begin{align}
		&\| \langle A\rangle^{-s} \varphi_{i,a}(\varepsilon A) G_{\pm\varepsilon}^p \varphi_{j,b}(\varepsilon A)
		\langle A\rangle^{-s} \| 
		\leq \|\langle A\rangle^{-s}\varphi_{i,a}(\varepsilon A)\| \|G_{\pm\varepsilon}\|^{p-1} \|G_{\pm\varepsilon}\langle A\rangle^{-s}\| \|\varphi_{j,b}\|_{L^\infty} \nonumber \\
		&\quad\leq C \varepsilon^s \varepsilon^{-p+1} \varepsilon^{-\frac12}
		\leq C \varepsilon^{s-\frac{1}{2}-\ell},
	\end{align}
	where we used that $\|\angles{A}^{-s} \varphi_{i,a}(\varepsilon A)\|\leq C \varepsilon^s$, which follows from $a+i\geq s$.
	
	It remains to estimate the sum over $n$ in \eqref{eq:IteratedRegularisationDerivative}. Expanding the commutators as before, we obtain the identity
	\begin{align} \label{eq:ExpansionCommutator}
		&\varphi_i(\varepsilon A) \mathcal A^{k-1-n} [ D_\zeta^\ell \partial_\varepsilon^n (\zeta G_{\pm\varepsilon} K_{\pm\varepsilon}^*G_{\pm\varepsilon}) ] \varphi_j(\varepsilon A) \nonumber \\
		&\quad = \sum_{\substack{a,b\geq 0 \\ a+b=k-1-n}} \frac{(k-1-n)!}{a!b!} (-1)^a \varepsilon^{-k+1+n} \varphi_{i,a}(\varepsilon A) D_\zeta^\ell \partial_\varepsilon^n (\zeta G_{\pm\varepsilon} K_{\pm\varepsilon}^*G_{\pm\varepsilon}) \varphi_{j,b}(\varepsilon A).
	\end{align}
	By the Leibniz formula, $D_\zeta^\ell \partial_\varepsilon^n (\zeta G_{\pm\varepsilon} K_{\pm\varepsilon}^*G_{\pm\varepsilon})$ is a linear combination of terms of the form
	\begin{align}
		\zeta (D_\zeta^p\partial_\varepsilon^u G_{\pm\varepsilon})
		(\partial_\varepsilon^wK_{\pm\varepsilon}^*)	(D_\zeta^q\partial_\varepsilon^v G_{\pm\varepsilon}),
	\end{align}
	where $u+v+w=n$ and $p+q\leq\ell$. For each such a term, we obtain the bound
	\begin{align}
		&\|\angles{A}^{-s} \zeta (D_\zeta^p\partial_\varepsilon^u G_{\pm\varepsilon})
		(\partial_\varepsilon^wK_{\pm\varepsilon}^*)	(D_\zeta^q\partial_\varepsilon^v G_{\pm\varepsilon}) \angles{A}^{-s}\| \nonumber \\
		&\quad \leq \| \angles{A}^{-s} D_\zeta^p\partial_\varepsilon^u G_{\pm\varepsilon} \| \| \partial_\varepsilon^wK_{\pm\varepsilon}^* \| \| D_\zeta^q\partial_\varepsilon^v G_{\pm\varepsilon} \angles{A}^{-s} \| \nonumber \\
		&\quad \leq C \varepsilon^{-p-u-\frac{1}{2}}
		\varepsilon^{s-\frac{1}{2}-w} \varepsilon^{-q-v-\frac{1}{2}} \nonumber \\
		&\quad \leq C \varepsilon^{s-\frac{3}{2}-p-q-n} \leq C \varepsilon^{s-\frac{3}{2}-\ell-n},
	\end{align}
	where we applied Proposition~\ref{prop:Regularity} and used the assumption $U\in \mathcal{C}^{s+1/2}(A)$ to estimate the derivatives of $K_{\pm\varepsilon}$, and Lemma~\ref{lem:DerivativesG} to estimate the derivatives of $G_{\pm\varepsilon}$. Inserting this estimate into \eqref{eq:ExpansionCommutator}, we obtain
	\begin{align}\label{eq:GlobalMixedErrorEstimate}
		\big\|\langle A\rangle^{-s} \varphi_i(\varepsilon A) \mathcal A^{k-1-n} [ D_\zeta^\ell \partial_\varepsilon^n (\zeta G_{\pm\varepsilon} K_{\pm\varepsilon}^*G_{\pm\varepsilon})] \varphi_j(\varepsilon A) \langle A\rangle^{-s} \big\| \leq C \varepsilon^{s-\frac{1}{2}-\ell-k}.	
	\end{align}
	This finishes the proof of \eqref{eq:HighestMixedDerivativeEstimate}. The estimate \eqref{eq:MixedDerivativeEstimate} is an immediate consequence of \eqref{eq:HighestMixedDerivativeEstimate}, by integrating \eqref{eq:HighestMixedDerivativeEstimate} with respect to $\varepsilon$ sufficiently many times.
\end{proof}

\section{Proofs of Main Results}
\label{sec:Proofs}

\begin{proof}[Proof of Theorem~\ref{thm:LAP}]
	The statement of the theorem is equivalent to the following: for all compact subsets $I_0 \subset \mu^A(U)$ and all $f,g\in\mathcal{H}_{1/2,1}$, the limits of $\scp{f}{(1-r\zeta U^*)^{-1}g}$ exist, as $r$ approaches 1 from below or above, uniformly in $\zeta \in I_0$ \cite[p.~300]{amrein1996}. Moreover, using the identity
	\begin{align}
		(1-r\zeta U^*)^{-1}(1+r\zeta U^*) = 2(1-r\zeta U^*)^{-1} - 1,
	\end{align}	
	we see that proving the existence of these limits is equivalent to establishing the existence of the limits of $\scp{f}{(1-r\zeta U^*)^{-1}(1+r\zeta U^*)g}$, as $r$ approaches 1 from below or above, uniformly in $\zeta\in I_0$.
	
	We now prove the existence of the latter limits. We may assume that $f=g$, and deduce the existence of the limits in the case $f\neq g$ from the polarisation method. We may also assume that $I_0$ is a sufficiently small compact set, satisfying the assumptions of Section~\ref{sec:Method}. Let $\varphi\in C_0^\infty(\mathbb{R})$ be a real-valued function which satisfies $\varphi=1$ in a neighbourhood of 0, and define $f_\varepsilon = \varphi(\varepsilon A) f$. Let $\zeta\in I_0$, $r\in (0,\infty)\backslash\{1\}$, and $\sigma$ be the sign of $1-r^2$. Then, by Proposition~\ref{prop:QuadraticEstimate}, there exists $\varepsilon_0>0$ such that $1-r\zeta(U_{\sigma\varepsilon})^{*}$ is invertible for all $\varepsilon\in[0,\varepsilon_0]$. As in Section~\ref{ssec:DifferentialInequality}, we define the function
	\begin{align}
		F_{\sigma\varepsilon}(r,\zeta) = \scp{f_\varepsilon}{(1-r\zeta(U_{\sigma\varepsilon})^*)^{-1}(1+r\zeta(U_{\sigma\varepsilon})^*)f_\varepsilon}.
	\end{align}
	Moreover, by Lemma~\ref{lem:RegularityFepsilon}, the function $\varepsilon \mapsto F_{\sigma\varepsilon}(r,\zeta)$ is of class $C^\infty$ on $(0,\varepsilon_0]$ and continuous on $[0,\varepsilon_0]$. Hence, by the fundamental theorem of calculus, we may write
	\begin{align}
		F_0(r,\zeta) = F_{\sigma\varepsilon_0}(r,\zeta) - \int_0^{\varepsilon_0} F_{\sigma\varepsilon}'(r,\zeta) \diff\varepsilon.
	\end{align}
	According to Corollary~\ref{cor:ExtensionsFepsilon}, the functions $F_{\sigma\varepsilon_0}(r,\zeta)$ have a limit as $r$ approaches 1 from below or above, uniformly in $\zeta\in I_0$. Moreover, by Proposition~\ref{prop:MajorisingFunction}, there exists a function $\rho\in L^1((0,\varepsilon_0))$ such that $\sup_{r\in [0,\infty), \zeta\in I_0} |F'_{\sigma\varepsilon}(r,\zeta)| \leq \rho(\varepsilon)$. According to Corollary~\ref{cor:ExtensionsFepsilon}, the functions $F'_{\sigma\varepsilon}(r,\zeta)$, for $\varepsilon\in (0,\varepsilon_0]$, have a limit as $r$ approaches 1 from below or above, uniformly in $\zeta\in I_0$. Thus, by the dominated convergence theorem, we can exchange the limit in $r$ with the integral over $\varepsilon$.
\end{proof}

\begin{proof}[Proof of Theorem~\ref{thm:HigherRegularity}]
	Let $\zeta_0\in \mu^A(U)$ and let $\mathbb{T} \neq I \subset \mu^A(U)$ be an open subset containing $\zeta_0$ as in Section~\ref{sec:Method}. Let $I_0\subset I$ be a compact subset. For $\zeta\in I_0$ and $\varepsilon\in (0,\varepsilon_0)$, we define, as in Section~\ref{ssec:HigherOrderEstimates}, the function
	\begin{align}
		L_\pm(\zeta,\varepsilon) = \varphi(\varepsilon A) G_{\pm\varepsilon}(1,\zeta) \varphi(\varepsilon A).
	\end{align}
	Let $m\geq 2s > 1$ be an integer. In Proposition~\ref{prop:MixedDerivativeEstimates}, we proved that, if $U\in\mathcal{C}^{s+1/2}(A)$, there exists a constant $C<\infty$ such that, for all $\varepsilon\in (0,\varepsilon_0)$, 
	\begin{align} 
		\sup_{\zeta\in I_0} \sum_{\ell+k=m} \| \langle A\rangle^{-s} D_\zeta^\ell\partial_\varepsilon^k	L_\pm(\zeta,\varepsilon) \langle A\rangle^{-s} \| \leq C \varepsilon^{s-\frac{1}{2}-m}.
	\end{align}
	From Proposition~\ref{prop:CircleBoundaryRegularity} in the appendix, it follows that the limits $b_\pm(\zeta) = \lim_{\varepsilon\downarrow0} L_\pm(\zeta,\varepsilon)$ exist in $B(\mathcal{H}_s,\mathcal{H}_{-s})$ and that they are locally of class $\Lambda^{s-1/2}$. By Theorem~\ref{thm:LAP} and its proof, we can identify the limits $b_\pm(\zeta)$ with $B_\pm(\zeta)$.
\end{proof}

\paragraph*{Acknowledgements.} I thank Serge Richard for his feedback on the content of the manuscript and for his hospitality during my stay at Nagoya University. I gratefully acknowledge financial support from the German Academic Exchange Service (DAAD) through the Postdoctoral Researchers International Mobility Experience (PRIME) programme. The proof of Proposition~\ref{prop:PhiInfty} was found with the help of ChatGPT.

\appendix

\section{Boundary Regularity Criterion}

Let $\mathcal{X}$ be a Banach space, $I\subset \mathbb{T}$ an open subset, and $\varepsilon_0>0$. We recall the boundary regularity criterion of \cite[Proposition~A.1]{boutetDeMonvel1996} adapted to the unit circle. 

\begin{prop}\label{prop:CircleBoundaryRegularity}
	Let $F:I\times(0,\varepsilon_0)\to \mathcal X$ be a function of class $C^m$ for some integer $m\geq1$, and let $\alpha\in(0,m)$. Assume that, for every compact subset $I_0\subset I$, there exists a constant $C_{I_0}$ such that, for all $\varepsilon\in (0,\varepsilon_0)$,
	\begin{align}\label{eq:CircleDerivativeHypothesis}
		\sup_{\zeta \in I_0} \sum_{\substack{\ell,k\geq 0 \\ \ell+k=m}} \|D_\zeta^\ell\partial_\varepsilon^k F(\zeta,\varepsilon)\|_{\mathcal X} \leq C_{I_0} \varepsilon^{\alpha-m}.
	\end{align}
	Then, the limit
	\begin{align}
		f(\zeta) = \lim_{\varepsilon\downarrow0} F(\zeta,\varepsilon)
	\end{align}
	exists in $\mathcal X$, locally uniformly in $\zeta$. Moreover, $f$ is locally of class $\Lambda^\alpha$.
\end{prop}

\begin{proof}
	We first prove the existence of the boundary values $f(\zeta)$. Let $I_0\subset I$ be a compact subset. From Taylor's formula, applied to the second variable, it follows that, for $\zeta\in I_0$ and $\varepsilon,\delta\in(0,\varepsilon_0)$ with $\varepsilon<\delta$,
	\begin{align}\label{eq:CircleTaylorFormula}
		F(\zeta,\varepsilon) = \sum_{k=0}^{m-1} \frac{(\varepsilon-\delta)^k}{k!} \partial_\varepsilon^kF(\zeta,\delta)
		+ \frac{(-1)^m}{(m-1)!} \int_{\varepsilon}^{\delta} (\tau-\varepsilon)^{m-1} \partial_\tau^m F(\zeta,\tau) \diff\tau,
	\end{align}
	where the integral is defined in the sense of Bochner. Furthermore, from the assumption \eqref{eq:CircleDerivativeHypothesis}, it follows that, for $\tau\in [\varepsilon,\delta]$,
	\begin{align}
		(\tau-\varepsilon)^{m-1} \sup_{\zeta\in I_0} \|\partial_\tau^mF(\zeta,\tau)\|_{\mathcal X} \leq C_{I_0} \tau^{m-1}\tau^{\alpha-m} = C_{I_0} \tau^{-1+\alpha}.
	\end{align}
	The function $\tau^{-1+\alpha}$ is integrable at zero because $\alpha>0$. Thus, by the dominated convergence theorem, we can take the limit $\varepsilon\downarrow0$ in \eqref{eq:CircleTaylorFormula}, and obtain the boundary value function
	\begin{align}\label{eq:CircleBoundaryRepresentation}
		f(\zeta) = \sum_{k=0}^{m-1} \frac{(-1)^k}{k!} \delta^k \partial_\varepsilon^kF(\zeta,\delta)
		+ \frac{(-1)^m}{(m-1)!}\int_0^\delta \tau^{m-1} \partial_\tau^m F(\zeta,\tau)\diff\tau.
	\end{align}	
	Clearly, the convergence is uniform in $\zeta\in I_0$. 
	
	We estimate the $m$-th finite differences of $f$. Let $I_1\subset I$ be a compact subset such that, for all~$\theta$ sufficiently close to 0, $\e^{\I m\theta} \zeta \in I_1$ for all $\zeta\in I_0$. Applying $\Delta_\theta^m$ to \eqref{eq:CircleBoundaryRepresentation} gives
	\begin{align}\label{eq:CircleDifferenceRepresentation}
		\Delta_\theta^m f(\zeta) = \sum_{k=0}^{m-1} \frac{(-1)^k}{k!} \delta^k \Delta_\theta^m \partial_\varepsilon^k F(\zeta,\delta)
		+ \frac{(-1)^m}{(m-1)!} \int_0^{\delta} \tau^{m-1} \Delta_\theta^m\partial_\tau^m F(\zeta,\tau) \diff\tau.
	\end{align}
	To estimate the integral term, we expand $\Delta_\theta^m$ using the definition \eqref{eq:FiniteDifference} and use \eqref{eq:CircleDerivativeHypothesis} for the compact set $I_1$. Then, for $\theta$ sufficiently close to 0,
	\begin{align}
		\sup_{\zeta\in I_0} \|\Delta_\theta^m \partial_\tau^m F(\zeta,\tau)\|_{\mathcal X} \leq C_{I_1} 2^m \tau^{\alpha-m}.
	\end{align}
	Integrating this bound, we obtain, for $\theta$ sufficiently close to 0,
	\begin{align}\label{eq:CircleIntegralDifferenceBound}
		\sup_{\zeta\in I_0} \left\| \frac{(-1)^m}{(m-1)!} \int_0^\delta \tau^{m-1} \Delta_\theta^m\partial_\tau^m F(\zeta,\tau) \diff\tau \right\|_{\mathcal X} 
		\leq C_{I_1} \frac{2^m}{(m-1)!\alpha}\delta^\alpha.
	\end{align}		
	Next, we estimate the sum over $k$ in \eqref{eq:CircleDifferenceRepresentation}. We use $\Delta_\theta^m=\Delta_\theta^k\Delta_\theta^{m-k}$, expand the finite difference of order $k$ according to its definition, and bound the finite difference of order $m-k$ by the tangential derivative of order $m-k$. This yields, for $\theta$ sufficiently close to 0,
	\begin{align}
		\sup_{\zeta\in I_0} \|\Delta_\theta^m \partial_\varepsilon^kF(\zeta,\delta)\|_{\mathcal X}
		\leq 2^k|\theta|^{m-k} \sup_{\zeta\in I_1} \|D_\zeta^{m-k}\partial_\varepsilon^k F(\zeta,\delta)\|_{\mathcal X}
		\leq C_{I_1} 2^k |\theta|^{m-k}\delta^{\alpha-m}.
	\end{align}		
	Combining these estimates, we find that, for $\theta$ sufficiently close to 0,
	\begin{align}
		\sup_{\zeta\in I_0} \|\Delta_\theta^m f(\zeta)\|_{\mathcal X} 
		\leq C_{I_1} \sum_{k=0}^{m-1} \frac{2^k}{k!} |\theta|^{m-k} \delta^{\alpha-m+k}
		+ C_{I_1} \frac{2^m}{(m-1)!\alpha}\delta^\alpha.
	\end{align}
	Finally, we choose $\delta=|\theta|$, which is permitted if $\theta$ is sufficiently close to 0. The resulting estimate proves that $f$ is locally of class $\Lambda^\alpha$.
\end{proof}

\section{Functions Preserving the Regularity Classes $C^k(A)$}
\label{appx:Counterexample}

For $k\in\mathbb{N}$, we denote by $\Phi^k$ the set of bounded Borel functions $\varphi$ such that $H\in C^k(A)$ implies $\varphi(H) \in C^k(A)$, for all pairs of self-adjoint operators $A,H$. It is well known that $C_0^\infty(\mathbb{R}) \subset \Phi^\infty = \bigcap_{k\in\mathbb{N}} \Phi^k \subset C^\infty(\mathbb{R})$ (see \cite[Theorem~6.2.5, Example~6.2.7]{amrein1996} for more precise results). 

In \cite[p.~248]{amrein1996}, it is asked whether a function $\varphi\in C^\infty(\mathbb{R})$ which satisfies $\varphi(x) = 0$ in a neighbourhood of $-\infty$ and $\varphi(x)=1$ in a neighbourhood of $\infty$ belongs to $\Phi^\infty$. We show that this is not the case by explicitly characterising the class $\Phi^\infty$. As in Section~\ref{ssec:SelfAdjointOperators}, the Cayley transform is denoted by $\kappa$. It extends to a diffeomorphism from $\mathbb{R} \cup \{\infty\}$ to $\mathbb{T}$ with $\kappa(\infty)=1$. We also recall that the regularity classes $C^k(A)$ for unbounded self-adjoint operators are understood in the sense that $H$ is of class $C^k(A)$ if its resolvent $(H+\I)^{-1}$ is of class $C^k(A)$.

\begin{prop} \label{prop:PhiInfty}
	$\Phi^\infty = \{\varphi\in C^\infty(\mathbb R): s\mapsto \varphi(1/s) \text{ extends smoothly through }s=0 \}.$
\end{prop}

\begin{proof}	
	Suppose that $\varphi\in\Phi^\infty \subset C^\infty(\mathbb{R})$. On $L^2(\mathbb{T})$, let $A$ be the generator of angular translations with Sobolev domain $H^1(\mathbb{T})$, and let $H=M_{\kappa^{-1}}$ be the multiplication operator by $\kappa^{-1}$ on its maximal domain. The resolvent of $H$ is given by $(H+\I)^{-1} = (1-M_{\zeta\mapsto\zeta})/(2\I)$. Furthermore, the identity
	\begin{align}
		\e^{\I \tau A} M_{\zeta\mapsto\zeta} \e^{-\I \tau A} = \e^{\I \tau} M_{\zeta\mapsto\zeta}, \qquad \tau\in\mathbb{R},
	\end{align}
	shows that $H\in C^\infty(A)$. Thus, by assumption on $\varphi$, the bounded operator $\varphi(H) = M_{\varphi \circ \kappa^{-1}}$ is also of class $C^\infty(A)$. Moreover, the constant function $1 \in L^2(\mathbb{T})$ is invariant under angular translations; hence,
	\begin{align}
		\e^{\I \tau A} \varphi(H) \e^{-\I \tau A} 1 = \e^{\I \tau A} M_{\varphi\circ \kappa^{-1}} \e^{-\I \tau A} 1 = \e^{\I \tau A} (\varphi\circ \kappa^{-1}).
	\end{align}
	This function is smooth in $\tau$, as $\varphi(H)\in C^\infty(A)$. Because $A$ is the self-adjoint generator of a unitary group, we deduce that $\varphi \circ \kappa^{-1} \in \bigcap_{m\ge1}D(A^m) = \bigcap_{m\ge1}H^m(\mathbb{T})$. By the Sobolev embedding theorem, $\varphi\circ \kappa^{-1}$ has a smooth representative $\psi\in C^\infty(\mathbb{T})$, which agrees with $\varphi\circ \kappa^{-1}$ on $\mathbb{T}\backslash \{1\}$. Writing $\varphi = \psi \circ \kappa$, we conclude that the map $s\mapsto \varphi(1/s)$ has a smooth extension through $s=0$. 
	
	The converse inclusion has already been proved in \cite[Corollary~3.6.13]{amrein1991}. Here, we give an independent proof, which uses the Cayley transform. Suppose that $\psi = \varphi \circ \kappa^{-1}$ extends to a smooth function on $\mathbb{T}$. Fix $k\in\mathbb{N}$ and a pair of self-adjoint operators $A,H$ with $H\in C^k(A)$. Let $U=\kappa(H) \in C^k(A)$ be the Cayley transform of $H$. Note that the inverse $U^{-1}$ is also of class $C^k(A)$. From the product rule, it follows that, for every integer $0\leq j\leq k$, there exists a constant $C_j$ such that, for all $n\in\mathbb{Z}$ and $\tau \in\mathbb{R}$,
	\begin{align} \label{eq:DerivativeEstimate}
		\|\partial_\tau^j \mathcal{W}_\tau [U^n]\| \leq C_j(1+|n|)^j,
	\end{align}
	where the group $\mathcal{W}_\tau$ was defined in Section~\ref{ssec:OperatorRegularity}. Moreover, using the Fourier expansion of $\psi$, we obtain the identity
	\begin{align} \label{eq:FourierExpansion}
		\mathcal{W}_\tau[\varphi(H)] = \mathcal{W}_\tau [\psi(U)] = \frac{1}{\sqrt{2\pi}} \sum_{n\in\mathbb Z} \widehat{\psi}(n) \mathcal{W}_\tau[U^n].
	\end{align}	
	Because $\psi$ is smooth on $\mathbb{T}$, the Fourier coefficients $\widehat{\psi}(n)$ decay rapidly in $n$. Together with \eqref{eq:DerivativeEstimate}, this implies that \eqref{eq:FourierExpansion} and its strong derivatives up to order $k$ converge in operator norm, uniformly in $\tau\in\mathbb{R}$. Thus, $\varphi(H)\in C^k(A)$.
\end{proof}

\end{document}